\documentclass[12 pt, onecolumn]{IEEEtran}
\usepackage{amsmath,amssymb,euscript,latexsym,graphicx}
\usepackage{bbm,color,amstext,wasysym,subfig,parskip,balance,hyperref,url}
\graphicspath{{./},{./figures/}}

\newtheorem{thm}{Theorem}
\newtheorem{cor}[thm]{Corollary}

\newtheorem{lemma}[thm]{Lemma}
\newtheorem{prop}{Proposition}

\newtheorem{remark}[thm]{Remark}

\newcommand{\ee}{\end{equation}}
\newcommand{\be}{\begin{equation}}

\newcommand{\mK}{\mathbb K}

\newcommand{\mC}{{\mathbb C}}

\newcommand{\cD}{{\mathcal D}}
\newcommand{\cH}{{\mathcal H}}

\newcommand{\cS}{{\mathcal S}}

\newcommand{\trace}{\operatorname{tr}}

\def\endproof{\ \hfill QED. \bigskip}
\def\proof#1{\noindent {\bf Proof{#1}:}}

\definecolor{grey}{rgb}{0.6,0.6,0.6}
\definecolor{lightgray}{rgb}{0.97,.99,0.99}

\def\spacingset#1{\def\baselinestretch{#1}\small\normalsize}
\spacingset{1.1}

\begin{document}
\title{On the Matrix Monge-Kantorovich Problem}

\author{Yongxin Chen, Wilfrid Gangbo, Tryphon T. Georgiou, and Allen Tannenbaum
\thanks{Y.\ Chen is with the Department of Medical Physics, Memorial Sloan Kettering Cancer Center, NY; email: chen2468@umn.edu}
\thanks{W.\ Gangbo is with the Department of Mathematics, UCLA, Los Angeles, CA; email: wgangbo@math.ucla.edu}
\thanks{T.\ T. Georgiou is with the Department of Mechanical and Aerospace Engineering, University of California, Irvine, CA; email: tryphon@uci.edu}
\thanks{A.\ Tannenbaum is with the Departments of Computer Science and Applied Mathematics \& Statistics, Stony Brook University, NY; email: allen.tannenbaum@stonybrook.edu}}

\maketitle

\begin{abstract}
The classical Monge-Kantorovich (MK) problem as originally posed is concerned with how best to move a pile of soil or rubble to an excavation or fill with the least amount
of work relative to some cost function. When the cost is given by the square of the Euclidean distance, one can define a metric on densities called the {\em Wasserstein distance}. In this note, we formulate a natural matrix counterpart of the MK problem for positive definite density matrices. We prove a number of results about this metric including showing that it can be formulated as a convex optimization problem, strong duality, an analogue of the Poincar\'e-Wirtinger inequality, and a Lax-Hopf-Oleinik type result.
\end{abstract}

\section{Introduction}

The mass transport problem was first formulated by Monge in 1781 and concerned
finding the optimal way, in the sense of minimal transportation cost, of moving a pile
of soil from one site to another. This problem of \textbf{\emph{optimal mass transport}} (OMT) was given a modern formulation in the
work of Kantorovich, and so is now known as the Monge--Kantorovich (MK) problem; see \cite{Rachev,Villani} and the many references therein. As originally formulated, the problem is  \textbf{\emph{static}}. Namely, given two probability densities, one can define a metric, now known as the {\em Wasserstein distance}, that quantifies the cost of transport and enjoys a number of remarkable properties as described in \cite{Rachev,Villani}.  Optimal mass transport is a very active area of research with applications to
numerous disciplines including probability, econometrics, fluid dynamics,
automatic control, transportation, statistical physics, shape
optimization, expert systems, and meteorology.

A major development in optimal mass transport theory was realized in the seminal dynamic approach to optimal mass transport by Benamou and Brenier \cite{French}. These authors base their approach to OMT on ideas from fluid mechanics via the minimization of a kinetic energy functional subject to a continuity constraint.
In a recent paper \cite{Yongxin} by Chen \emph{et al}., a non-commutative counterpart of optimal transport was developed where density matrices $\rho$ (i.e., Hermitian matrices that are positive-definite and have unit trace) replace probability distributions, and where ``transport'' corresponds to a flow on the space of such matrices that minimizes a corresponding action integral, thereby extending the fluid dynamical approach of \cite{French}. (A similar approach to \cite{Yongxin} was done independently by Carlen and Maas at about the same time \cite{Carlen}.) Here again, based on the continuity equation that imposes a ``mass preservation'' constraint on a quadratic ``kinetic energy,'' we study a convex optimization problem that leads to a certain Riemannian structure on densities matrices and generalizes the work of \cite{Otto} to the current non-commutative setting.

The present paper is a continuation of \cite{Yongxin} in which a number of the results are given rigorous mathematical proofs. We show that indeed in line with \cite{French}, one has a convex optimization problem and strong duality, an analogue of the Poincar\'e-Wirtinger inequality, and a Lax-Hopf-Oleinik type result, all in our non-commutative Wasserstein framework.

\section{Continuity equation \& Wasserstein distance} \label{sec:quantum}

In this section, we set up the continuity equation that is the basis for our formulation of the Wasserstein distance for density matrices. We follow closely the recent paper \cite{Yongxin}. In that work, an approach is developed based on the \textbf{\emph{Lindblad equation}} which describes the evolution of open quantum systems. Open quantum systems are thought of as being coupled to a larger system (heat bath), and thus cannot in general be described by a wave
function and a unitary evolution. The proper description is in terms of a density operator $\rho$ \cite{Sigal} which in turn obeys the Lindblad equation where we assume $\hbar =1$:
	\begin{equation}\label{eq:lindblad}
	\dot\rho =-i [H, \rho]+\sum_{k=1}^N( L_k\rho L_k^*-\frac{1}{2}\rho L_k^*L_k-\frac{1}{2}L_k^*L_k\rho).
	\end{equation}
Here, * as superscript denotes conjugate transpose and $[H, \rho] := H \rho-\rho H$ denotes the commutator.
The first term on the right-hand side describes the evolution of the state under the effect of the Hamiltonian $H$, and it is unitary (energy preserving), while the other the terms on the right-hand side model diffusion and, thereby, capture the dissipation of energy; these dissipative terms together represent the quantum analogue of Laplace's
operator $\Delta$  (as it will become clear shortly). The Lindblad equation defines a non-unitary evolution of the density and the calculus we develop next actually underscore parallels with classical diffusion and the Fokker-Planck equation.

Denote by $\cH$ and $\cS$ the set of  $n\times n$ Hermitian and skew-Hermitian matrices, respectively. We will assume that all of our matrices are fixed to be $n\times n$. Next, we denote the space of block-column vectors consisting of $N$ elements in $\cS$ and $\cH$ as $\cS^N$, respectively $\cH^N$. We let $\cH_+$ and $\cH_{++}$ denote the cones of nonnegative and positive-definite matrices, respectively, and
	\begin{eqnarray}\label{eq:D}
    \cD &:=&\{\rho \in \cH_{+} \mid \trace(\rho)=1\},\\
	\cD_+ &:=&\{\rho \in \cH_{++} \mid \trace(\rho)=1\}.
    \end{eqnarray}
We note that the tangent space of $\cD_+,$ at any $\rho\in \cD+$ is given by
	\begin{equation}\label{eq:Trho}
		T_\rho=\{ \delta \in \cH \mid \trace(\delta)=0\},
	\end{equation}
and we use the standard notion of inner product, namely
	\[
		\langle X;Y\rangle=\trace(X^*Y),
	\]
for both $\cH$ and $\cS$.
For $X, Y\in \cH^N$ ($\cS^N$),
	\[
		\langle X; Y\rangle=\sum_{k=1}^N \trace(X_k^*Y_k).
	\]
Given $X=[X_1^*,\cdots,X_N^*]^* \in \cH^N$ ($\cS^N$), $Y\in \cH$ ($\cS$), set
	\[
		XY=\left[\begin{array}{c}
		X_1\\
		\vdots \\
		X_N
		\end{array}\right]Y
		:=
		\left[\begin{array}{c}
		X_1Y\\
		\vdots \\
		X_N Y
		\end{array}\right],
	\]
and
	\[
		YX=Y\left[\begin{array}{c}
		X_1\\
		\vdots \\
		X_N
		\end{array}\right]
		:=
		\left[\begin{array}{c}
		YX_1\\
		\vdots \\
		YX_N
		\end{array}\right].
	\]

If we assume that $L_k=L_k^*$, i.e., $L_k\in \cH$ for all $k\in 1 \ldots, N$, then we can define
	\begin{equation}\label{eq:gradient}
		\nabla_L: \cH \rightarrow {\cS}^N, ~~X \mapsto
		\left[ \begin{array}{c}
		L_1 X-XL_1\\
		\vdots \\
		L_N X-X L_N
		\end{array}\right]
	\end{equation}
as the \textbf{\emph{gradient operator}}.
Note that $\nabla_L$ acts just like the standard gradient operator, and in particular, satisfies
	\be \label{eq:product}
		\nabla_L(XY+YX)=(\nabla_L X) Y+ X (\nabla_L Y)+(\nabla_L Y) X+ Y (\nabla_L X),~~\forall X, Y\in \cH.
    \ee
The dual of $\nabla_L$, which is an analogue of the (negative) \textbf{\emph{divergence operator}}, is defined by
	\begin{equation}\label{eq:divergence}
		\nabla_L^*: {\cS}^N \rightarrow \cH,~~Y=
		\left[ \begin{array}{c}
		Y_1\\
		\vdots \\
		Y_N
		\end{array}\right]
		\mapsto
		\sum_k^N L_k Y_k-Y_k L_k.
	\end{equation}
The duality
	\[
		\langle \nabla_L X ; Y\rangle =\langle X ; \nabla_L^* Y\rangle
	\]
follows by definition.

With these definitions, we define the (matricial) \textbf{\emph{Laplacian}} as
	\[
		\Delta_L X:=-\nabla_L^*\nabla_L X=\sum_{k=1}^N( 2L_k X L_k^*-X L_k^*L_k-L_k^*L_k X),~~X\in \cH,
	\]
which is exactly (after scaling by $1/2$) the diffusion term in the Lindblad equation  \eqref{eq:lindblad}. Therefore Lindblad's equation (under the assumption that $L_k=L_k^*$) can be rewritten as
	\[
		\dot\rho =-i [H, \rho]+\frac{1}{2} \Delta_L \rho,
	\]
i.e., as a continuity equation expressing flow under the influence of a suitable vector field.

In our case, we will consider a \textbf{\emph{continuity equation}} of the form
\[
\dot\rho =\nabla^*_L m
\]
for $m\in \cS^N$ a suitable ``momentum field.'' In particular, we are interested in
the following family of continuity equations:
\begin{equation}\label{eq:continuitygen}
	\dot \rho=\nabla_L^*  M_\rho(v), 
	\end{equation}
where the momentum field is expressed as a \textbf{\emph{non-commutative product}}
$M_\rho(v)\in \cS^N$ between a ``velocity field'' $v\in \cS^N$ and the density matrix $\rho$.

Several such ``non-commutative products'' have been considered (see \cite{Yongxin}), however, in the present work, we consider the following case:
\begin{subequations}
\begin{equation}\label{eq:product1}
		M_\rho(v):=\frac{1}{2}(v \rho+\rho v),
\end{equation}
which gives
\begin{equation}\label{eq:continuity}
		\dot \rho=\frac{1}{2} \nabla_L^* (v \rho+\rho v)
	\end{equation}
\end{subequations}
and $v=[v_1^*,\ldots,v_N^*]^* \in \cS^N$. Clearly $v \rho+\rho v \in \cS^N$, which is consistent with the definition of $\nabla_L^*$. In \cite{Yongxin}, we call this the \emph{anti-commutator} case, since
\[
v \rho+\rho v=:\{v,\rho\}
\]
is the anti-commutator when applied to elements of an associative algebra. In \cite{Yongxin}, another possibility is considered for the multiplication operator $M_\rho(v).$


Given two density matrices $\rho_0, \rho_1 \in \cD_+$ we formulate the optimization problem (following \cite{Yongxin})
 	\begin{subequations}\label{eq:quantumomt}
 	\begin{eqnarray}\label{eq:quantumomt1}
	 W_{2}(\rho_0, \rho_1)^2:=&&\inf_{\rho\in \cD_+, v\in \cS^N} \int_0^1 \trace(\rho v^*v) dt,\\
	&&\dot \rho=\frac{1}{2} \nabla_L^* (v \rho+\rho v), \label{eq:quantumomt2}\\
	&& \rho(0)=\rho_0, ~~\rho(1)=\rho_1,  \label{eq:quantumomt3}
	\end{eqnarray}
	\end{subequations}
and define the \textbf{\emph{Wasserstein distance}} $W_{2}(\rho_0, \rho_1)$ between $\rho_0$ and $\rho_1$ to be the square root of the infimum of the cost \eqref{eq:quantumomt1}.
Other choices for $M_\rho(v)$ in \eqref{eq:product1} give alternative Wasserstein metrics, as noted in \cite{Yongxin}.
In order for the metric $W_2$ to be well-defined for all $\rho_0, \rho_1 \in \cD_+$, we need to assume that ${\rm ker}(\nabla_L)$ is spanned by the identity matrix. The results in the present paper, however, carry through without this assumption.


\section{Quadratic forms and Poincar\'e-Wirtinger inequality} \label{sec:mk1}

In this section, we prove some initial convexity results as well as a Poincar\'e-Wirtinger type inequality that we will need in the sequel. We begin with some notation.
If $m_1, \cdots, m_N \in \mathbb C^{n \times n}$, we define the column vector $m\in \mathbb C^{nN \times n}$
with matrix entries the $m_i$'s
by
\[
m=(m_1^*, \cdots, m_N^*)^*,
\]
the column vector $m_*\in \mathbb C^{nN \times n}$ with entries $m_i^*$'s, by
\[\quad m_* =(m_1, \cdots, m_N)^*.
\]
For $m,b\in \mathbb C^{nN \times n}$, i.e., with
$b=(b_1^*, \cdots, b_N^*)^*$ for $b_1, \cdots, b_N \in \mathbb C^{n \times n}$, we define the inner products
\[
\langle m_i; b_i \rangle=  \trace (m^*_i b_i) , \quad  \langle m; b \rangle=  \trace(m^* b)= \sum_{i=1}^N \langle m_i; b_i \rangle
\]
and introduce
\[
 m \cdot b= {1 \over 2} \bigl(\langle m; b \rangle+\langle b; m \rangle \bigr) \in \mathbb R.
\]
Then, for $v \in \mathbb C^{nN \times n}$ and $\rho \in \cH_+$, we define the quadratic form
\begin{equation}\label{eq:quadraticform}
Q_\rho(v) := \trace(\rho v^* v)= \langle v \rho; v \rangle.
\end{equation}
The following lemma is an easy consequence of the definition and can be readily verified.

\begin{lemma}\label{rem:nov18.2016.0} Let $\rho \in \cH_+$  and $v, w \in \mathbb C^{nN \times n}$. The following hold:
\begin{enumerate}
\item[(i)] $Q_\rho(v) \geq 0$ and $Q_\rho(v) = 0$ if $v=0$; when $\rho\in\cH_{++}$, it becomes if and only if,
\item[(ii)]
$Q_\rho(v+w)=Q_\rho(v)+Q_\rho(w) + \langle v \rho;w \rangle +\langle w \rho; v \rangle$,
\item[(iii)]
$
Q_\rho \bigl((1-t) v + t w \bigr) =(1-t) Q_\rho(v) + t  Q_\rho(w) - t(1-t) Q_\rho (v-w),
$
\item[(iv)]
If we further  assume that $v, w \in \mathcal S^N$, then $\langle w \rho; v \rangle=\langle \rho v; w \rangle$.
\end{enumerate}
\end{lemma}

Since $|| \cdot||$ (the standard norm on $\cH$) is uniformly convex and ${\rm ker}(\nabla_L)$ is a finite dimensional space, there exists a unique ${\rm proj}(X) \in  {\rm ker}(\nabla_L)$, the orthogonal projection of $X$ onto ${\rm ker}(\nabla_L),$ such that
\[\min_{Z \in {\rm ker}(\nabla_L)} ||X-Z||=||X-{\rm proj}(X)||.\]
If we denote by ${\rm ker}(\nabla_L)^\perp$ the orthogonal complement of ${\rm ker}(\nabla_L)$ in $\mathcal H$, then
\[\mathcal H= {\rm ker}(\nabla_L) \oplus {\rm ker}(\nabla_L)^\perp.\]
%

\begin{lemma}\label{lem:nov18.2016.2bis} For any $\rho \in \mathcal H_+$, the map $X \rightarrow Q_\rho(\nabla_L X)$ is convex on $\mathcal H$. If in addition $\rho>0,$ then the map is strictly convex on ${\rm ker}(\nabla_L)^\perp$.
\end{lemma}
\proof{} From Lemma~\ref{rem:nov18.2016.0}, we obtain that, for $t \in (0,1)$ and $X, Y \in \mathcal H$,
\[
Q_\rho \bigl((1-t) \nabla_L X+ t \nabla_L Y \bigr)= (1-t)Q_\rho(\nabla_L X)+t Q_\rho(\nabla_L Y) - t(1-t) Q_\rho(\nabla_L X-\nabla_L Y).
\]
The convexity follows from the fact that $Q_\rho(\cdot)\ge 0$. Furthermore, if $\rho>0$, then $Q_\rho(\nabla_L X-\nabla_L Y)>0$ unless $\nabla_L X-\nabla_L Y=0.$ Hence,  $X \rightarrow Q_\rho(\nabla_L X)$ strictly convex on ${\rm ker}(\nabla_L)^\perp$. \endproof

\begin{thm}[Poincar\'e--Wirtinger inequality]\label{thm:nov18.2016.3} Let $\mathbb K \subset \cD_+$ be a compact set. Then there exists a constant $c_{\mathbb K}>0$ such that for all $X \in \mathcal H$ and $\rho \in \mathbb K$,
\[
Q_\rho \bigl(\nabla_L (X-{\rm proj}(X)  )\bigr)   \geq c_{\mathbb K} \bigl\| X-{\rm proj}(X)\bigr\|^2.
\]
\end{thm}
\proof{} Define
\begin{equation}\label{eq:nov18.2016.4b}
 c_{\mathbb K}:= \inf_{\rho, X} \Bigl\{ \trace \bigl(\rho (\nabla_L X)^* \nabla_L X \bigr) \; | \; \rho \in \mathbb K, X \in {\rm ker}(\nabla_L)^\perp,  ||X||=1\Bigr\}
 \end{equation}
and let $(\rho_k, X_k)_k$ be a minimizing sequence in (\ref{eq:nov18.2016.4b}). This infimum of a continuous function over a compact set is a minimum, attained at a certain $(\rho, X).$ Since   $ X \in {\rm ker}(\nabla_L)^\perp$, we cannot have $ X \in {\rm ker}(\nabla_L)$ because, we would otherwise have $X=0$ which will contradict the fact that $||X||=1.$ Since $\rho>0$ we conclude that
\begin{equation}\label{eq:nov18.2016.4}
\trace \Bigl(\rho (\nabla_L X)^* \nabla_L X \Bigr)>0
\end{equation}
In conclusion
 \[
 \trace \bigl(\rho (\nabla_L Y)^* \nabla_L Y \bigr) \geq  c_{\mathbb K}>0
 \]
 for any $\rho \in \mathbb K$ and any $ Y \in {\rm ker}(\nabla_L)^\perp$ such that  $||Y||=1.$  By homogeneity, this completes the proof of the theorem. \endproof

\begin{lemma}\label{le:dec20.2016.5} Given $\rho \in \cH_{++}$ and $v \in \mathcal S^N,$ there is a unique element $\nabla_L X, X\in\cH$,  such that
\[
Q_\rho(v-\nabla_L X) \leq Q_\rho(v-\nabla_L Y)
\]
for all $Y \in \cH.$ Furthermore,
the minimizer is characterized by the Euler--Lagrange equations
\begin{equation}\label{eq:dec20.2016.1}
(v-\nabla_L X) \rho + \rho(v-\nabla_L X)  \in {\rm ker} (\nabla^*_L).
\end{equation}

\end{lemma}
\proof{} Let $(X_\ell)_\ell \subset \cH$ be a sequence such that
\[
\lim_{\ell \rightarrow \infty} Q_\rho(v-\nabla_L X_\ell)= \inf_{Y \in \cH} Q_\rho(v-\nabla_L Y).
\]
Note that $\bigl(Q_\rho(\nabla_L X_\ell) \bigr)_\ell$ is bounded by definition and the convexity of $Q_\rho(\cdot)$. Replacing $X_\ell$ by $X_\ell- {\rm proj}(X_\ell) $ if necessary, we use the Poincar\'e--Wirtinger inequality to conclude that $(X_\ell)_\ell$ is a bounded sequence. Passing to a subsequence if necessary, we may assume that $(X_\ell)_\ell$ converges to some $X$ which minimizes  $Q_\rho(v-\nabla_L X)$ over $\cH.$  The uniqueness follows from the strict convexity of $Q_\rho(\cdot)$, and condition \eqref{eq:dec20.2016.1} expresses stationarity.
\endproof

\section{Flow rates in the space of densities}

We now return to the continuity equation
\[
\dot\rho =f
\]
with the {\em flow rate} $f$ being the divergence of a momentum field $p$, i.e.,
\begin{equation}\label{eq:flowrate}
f=\nabla^*_L\, p,
\end{equation}
with  $p\in \cS^N$, so that $f\in\cH$ as well as $\trace(f)=0$.
In particular, we are interested in the case where the momentum is a linear function of $\rho$ of the form $p=M_\rho(v)$ (see \eqref{eq:product1}); then
$p=\frac12(m-m_*)$ with $m = v\rho\in \mC^{Nn\times n}$, $\rho \in \cH_{++}$, and
$v \in \cS^N$, and 
\[
f = \frac12\nabla^*_L(v \rho + \rho v).
\]
Since the range of $\nabla^*_L$ coincides with ${\rm ker}(\nabla_L)^\perp$,  any $f$ belongs to ${\rm ker}(\nabla_L)^\perp$.
The next theorem states that in this case not only the converse holds, namely, that given $\rho \in \cH_{++}$, any $f \in {\rm ker}(\nabla_L)^\perp$
can be written as above with $m = v\rho,$ but that $v \in \cS^N$ can be selected in the range of $\nabla_L$ and that this
choice is unique.
%

\begin{thm} \label{thm:nov18.2016.6} For any $\rho \in \mathcal H_{++}$ and  $f \in {\rm ker}(\nabla_L)^\perp,$ there exists a unique $X \in {\rm ker}(\nabla_L)^\perp$  such that
\begin{equation}\label{eq:nov18.2016.4c}
f=\frac12\nabla^*_L\bigl(\nabla_L X\rho + \rho\nabla_L X \bigr).
 \end{equation}
 Furthermore, if $\mathbb K$ is a compact subset of $\cD_+$ and $\rho / \trace(\rho) \in \mathbb K$ then there exists $c_\mK>0$ such that
 \begin{equation}\label{eq:nov18.2016.4d}
 ||f|| \geq c_{\mathbb K} \; \trace(\rho) \, ||X||.
 \end{equation}
\end{thm}
\proof{} Define  the functional
\[
I(Y)=\frac12 Q_\rho(\nabla_L Y)- \langle f; Y \rangle, \quad \forall ~ Y \in  \cH.
\]
To avoid trivialities, we assume that $f \not =0.$ Observe that $I(Y) \equiv 0$ on ${\rm ker}(\nabla_L)$ and for $0<\epsilon<<1$ we have that
\[
I(\epsilon f)= \frac{\epsilon^2}{2} Q_\rho(\nabla_L f) -\epsilon ||f||^2<0
\]
 Thus, $\lambda_0$, the infimum of $I$ over $\mathcal H$ is negative. Let $(Y_\ell)_\ell$ be a minimizing sequence. Since $\lambda_0<0$, for $\ell$ large enough, $I(Y_\ell)<0$ and so, $Y_\ell \in \mathcal H \setminus {\rm ker}(\nabla_L).$ Replacing $Y_\ell$ by $Y_\ell-{\rm proj}(Y_\ell)$ if necessary, we may assume that $Y_\ell \in {\rm ker}(\nabla_L)^\perp.$ By the Poincar\'e--Wirtinger inequality
 \[
 0 > I(Y_\ell) \geq c_\mathbb K \trace(\rho) ||Y_\ell||^2-||f|| \, ||Y_\ell||.
 \]
Consequently, $(Y_\ell)_\ell$ is a bounded sequence and so, passing to a subsequence if necessary, we may assume that $(Y_\ell)_\ell$ converges to some $X \in {\rm ker}(\nabla_L)^\perp.$ We have $0 \geq c_\mathbb K \trace(\rho) ||X||^2-||f|| \, ||X||$ and so, (\ref{eq:nov18.2016.4d}) holds.

If $Y \in \mathcal H$ is arbitrary, then for any real number $\epsilon$, we use  Lemma \ref{rem:nov18.2016.0}
to conclude that
\[
I(X+ \epsilon Y)=I(X) +\epsilon  \bigl\langle \frac12\nabla_L^* \bigl(\nabla_L X \rho+\rho \nabla_L X \bigr)-f ; Y \bigr\rangle +o(\epsilon).
\]
By Lemma \ref{lem:nov18.2016.2bis}, $I$ is  convex on $\mathcal H$ and so, $X$ is a critical point of $I$ if and only if $X$ minimizes $I.$ Thus (\ref{eq:nov18.2016.4c}) holds if and only if $X$ minimizes $I.$ Since, the same lemma gives that $I$ is  strictly convex on ${\rm ker}(\nabla_L)^\perp$, $I$ admits a unique minimizer on ${\rm ker}(\nabla_L)^\perp$ which means that there exists a unique $X \in {\rm ker}(\nabla_L)^\perp$ such that (\ref{eq:nov18.2016.4c}) holds.  \endproof

\begin{remark} The uniqueness and existence of the representation may also be proven as follows. We first note that provided $\rho \in \cH_{++}$, the non-commutative multiplication in (\ref{eq:product1}), defines a positive definite
Hermitian operator
$$M_\rho : \cS^N \rightarrow \cS^N \; : \; v \mapsto \frac12(v \rho+\rho v).$$
It follows that
$\nabla^*_L M_\rho \nabla_L,$ when restricted to ${\rm ker}(\nabla_L)^\perp = {\rm range}(\nabla^*_L),$
is positive and therefore invertible. Thus, for all
$f \in {\rm ker}(\nabla_L)^\perp,$ (\ref{eq:nov18.2016.4c}) has a unique solution $X \in {\rm ker}(\nabla_L)^\perp$. One also gets a lower bound on the norm of $f$ as follows. If $\lambda_{min}>0$ denotes the smallest eigenvalue of
$$ \nabla^*_L M_\rho \nabla_L|_{{\rm ker}(\nabla_L)^\perp}: {\rm ker}(\nabla_L)^\perp \rightarrow {\rm ker}(\nabla_L)^\perp,$$
then
$$\|f \| \ge \lambda_{\rm min} \|X\|.$$
\end{remark}

\section{Flows in the space of densities}

We begin with establishing a canonical representation of flow rates that minimize a certain analogue of {\em kinetic energy} of our matrix-valued flows.

\begin{prop} \label{prop:dec18.2016.1} Suppose $\rho \in \mathcal H_{+}$,  $f \in {\rm ker}(\nabla_L)^\perp$, $X \in \cH$ satisfy  (\ref{eq:nov18.2016.4c}), and that $v \in \mathcal S^N$ is such that
\[
f={ 1 \over 2}\nabla^*_L\bigl(v \rho+\rho v \bigr).
\]
The following hold:
\begin{enumerate}
\item[(i)] For all $Y\in\cH$,
\[
{1 \over 2} \trace (\rho v^* v) \geq \langle f; Y \rangle-{1 \over 2} Q_\rho(\nabla_L Y)
\]
and equality holds if and only if $v  = \nabla_L X $.
\item[(ii)] Further assume that $\rho>0$ (which by Theorem \ref{thm:nov18.2016.6} is a sufficient condition for (\ref{eq:nov18.2016.4c}) to hold). Then
\[
\min_{m \in \mathbb C^{nN \times n}} \Bigl\{{1 \over 2} \langle m; m\rho^{-1} \rangle\; | \; f={1 \over 2}\nabla^*_L(m -m_*)  \Bigr\}= \max_{Y \in \cH} \left\{\langle f; Y \rangle-{1 \over 2} Q_\rho(\nabla_L Y)\right\}.
\]
Besides, the maximum is uniquely attained by the $X$ which satisfies (\ref{eq:nov18.2016.4c}) and the minimum is uniquely attained by $m= \nabla_L X \rho$ for the same $X$.
\end{enumerate}
\end{prop}
\proof{} (i) We have
\[
{1 \over 2} \trace (\rho v^* v)= {1 \over 2} \langle v\rho^{1 \over 2}; v \rho^{1 \over 2} \rangle + \langle f-{ 1 \over 2} \nabla^*_L\bigl(v \rho+\rho v \bigr) ; Y \rangle=
{1 \over 2} \langle v\rho^{1 \over 2};  v\rho^{1 \over 2} \rangle + \langle f; Y \rangle- { 1 \over 2}\langle  v \rho+\rho v ; \nabla_L Y \rangle.
\]
Since both $v$ as well as $\nabla_LY$ belong to $\cS^N$ (cf.\ Lemma \ref{rem:nov18.2016.0} (iv)), 
\[
\langle  v \rho+\rho v ; \nabla_L Y \rangle= \langle v  \rho^{1 \over 2} ; \nabla_L Y  \rho^{1 \over 2} \rangle + \langle \nabla_L Y  \rho^{1 \over 2}  ; v  \rho^{1 \over 2} \rangle.
\]
We conclude that
\[
{1 \over 2} \trace (\rho v^* v)= {1 \over 2} \| v \rho^{1 \over 2}- \nabla_L Y  \rho^{1 \over 2} \|^2  + \langle f; Y \rangle  -{1 \over 2} \|\nabla_L Y  \rho^{1 \over 2}\|^2 \geq  \langle f; Y \rangle  -{1 \over 2} Q_\rho(\nabla_L Y).
\]
(ii) Computations similar to the ones in (i) reveal that
\begin{equation}\label{eq:12.24.2016.1}
{1 \over 2} \langle m; m\rho^{-1} \rangle= {\frac{1}{2}} \| m \rho^{-\frac{1}{2}}- \nabla_L Y  \rho^{\frac{1}{2}}  \|^2  + \langle f; Y \rangle  -{\frac{1}{2}} \|\nabla_L Y  \rho^{ \frac{1}{2}}\|^2
\end{equation}
and so, for $Y\in\cH$,
\[
{1 \over 2} \langle m; m\rho^{-1} \rangle \geq  \langle f; Y \rangle  -{1 \over 2} Q_\rho(\nabla_L Y).
\]
 \endproof

We proceed to consider
 paths $\rho(t)\in \cH_+$ for $t\in[0,1]$ along with corresponding flow rates and action integrals. A corollary of the above proposition ascertains the measurability of the canonical representation of the velocity field $v$.

\begin{cor} \label{cor:dec18.2016.1}  Let $\mathbb L \subset \mathcal H_{++}$ and denote by  $A: {\rm ker}(\nabla_L)^\perp \times \mathbb L \rightarrow \cH$ the map  which to $(f, \rho)$ associates $X \in {\rm ker}(\nabla_L)^\perp$ such that  (\ref{eq:nov18.2016.4c}) holds.
\begin{enumerate}
\item[(i)] If $\mathbb L$ is a compact subset of $\mathcal H_{++}$, then $A$ is continuous.
\item[(ii)] If $\rho:[0,1] \rightarrow \cH_{++}$ and $f:[0,1] \rightarrow  {\rm ker}(\nabla_L)^\perp$ are continuous at $t_0 \in [0,1]$ then $A(f, \rho)$ is continuous at $t_0.$
\item[(iii)] If $\rho \in L^1(0,1; \cH_{++})$ and $f \in L^1(0,1;   {\rm ker}(\nabla_L)^\perp)$ are measurable, then  $A(f, \rho)$ is measurable.
\end{enumerate}

\end{cor}
\proof{} (i) Let $\mathbb K$ be the set of $\rho/ \trace(\rho)$ such that $\rho \in \mathbb L.$ Let $(f_\ell, \rho_\ell)_\ell$ be a sequence in ${\rm ker}(\nabla_L)^\perp \times \mathbb L$ converging to $(f, \rho).$ By Theorem \ref{thm:nov18.2016.6}, $(X_\ell)_\ell:=\bigl(A(f_\ell, \rho_\ell) \bigr)_l$ is a bounded sequence in  ${\rm ker}(\nabla_L)^\perp$ and so, has all its points of accumulation in ${\rm ker}(\nabla_L)^\perp.$ If $X$ is any such point of accumulation, then clearly
\[
f={ 1 \over 2}
\nabla_L^* (\nabla_L X \rho+ \rho \nabla_L X).
\]
Since $\rho$ is invertible, $X$ is unique and so, $A(f, \rho)=X.$ This establishes (i).

(ii) Condition (ii) is a direct consequence of (i).

(iii) Approximate $\rho$ in the $L^1$--norm by a sequence $(\rho_\ell)_\ell \subset C([0,1]; \cH_{++})$ which converges pointwise almost everywhere to $\rho.$ Similarly, approximate $f$ in the $L^1$--norm by a sequence $(f_\ell)_\ell \subset C([0,1]; {\rm ker}(\nabla_L)^\perp )$ which converges pointwise almost everywhere to $f.$ By (i),  $\bigl(A(f_\ell, \rho_\ell) \bigr)_\ell$ converges pointwise almost everywhere to $A(f, \rho)$ and so, $A(f, \rho)$ is measurable. This establishes (iii).
\endproof

\begin{lemma}\label{lem:nov18.2016.7}
If $\rho_0, \rho_1 \in \cD_+$, then the following hold:
\begin{enumerate}
\item[(i)] If $\rho_1 -\rho_0 \in   {\rm ker}(\nabla_L)^\perp$, then there exists a Borel map $t \rightarrow X(t) \in  {\rm ker}(\nabla_L)^\perp$ and a Borel map $t \rightarrow \rho(t)$ starting at $\rho(0)=\rho_0$ and ending at $\rho(1)=\rho_1$ such that
\begin{equation}\label{eq:nov19.2016.1}
\dot \rho(t) =\frac12\nabla^*_L\bigl(\nabla_L X(t) \rho(t)+\rho(t) \nabla_L X(t) \bigr)
\end{equation}
in the sense of distributions, and
\begin{equation}\label{eq:nov19.2016.1b}
\int_0^1 Q_{\rho(t)}( \nabla_L X(t)) dt < \infty.
\end{equation}
\item[(ii)] Conversely, assume that there exist a Borel map $t \rightarrow v(t) \in  \mathcal S^N$ and a Borel map $t \rightarrow \rho(t)$ starting at $\rho_0$ and ending at $\rho_1$ such that
\begin{equation}\label{eq:nov19.2016.2}
\dot \rho(t) =\frac12\nabla^*_L\bigl(v(t) \rho(t)+ \rho(t) v(t)\bigr)
\end{equation}
in the sense of distributions, 
and
\begin{equation}\label{eq:nov19.2016.2b}
\int_0^1 \trace( \rho(t) v(t)^*v(t)) dt < \infty.
\end{equation}
Then $\rho_1 -\rho_0 \in   {\rm ker}(\nabla_L)^\perp$.
\end{enumerate}
\end{lemma}
\proof{} (i) Assume $\rho_1 -\rho_0 \in   {\rm ker}(\nabla_L)^\perp$, Set $\rho(t)=(1-t) \rho_0+t \rho_1.$ Then $\mathbb K:=\{\rho(t) \; | \; t\in [0,1] \}$ is a compact subset of $\cD_+.$ For each $t \in[0,1]$, we use Theorem \ref{thm:nov18.2016.6}  to find a unique $X(t) \in {\rm ker}(\nabla_L)^\perp$ such that
\[
\rho_1 -\rho_0 =\frac12\nabla^*_L\bigl(\nabla_L X(t) \rho(t)+\rho(t) \nabla_L X(t)\bigr)
\]
and
\[
||X(t)\| \leq c_{\mathbb K}.
\]
By Corollary \ref{cor:dec18.2016.1}, $t \rightarrow X(t) \in  {\rm ker}(\nabla_L)^\perp$ is continuous. Hence, (\ref{eq:nov19.2016.1}) and (\ref{eq:nov19.2016.1b}) hold.

(ii) Conversely, assume (\ref{eq:nov19.2016.2}) and (\ref{eq:nov19.2016.2b}) hold. Let $Y \in {\rm ker}(\nabla_L)$, then
\[
\langle \rho_1 -\rho_0; Y \rangle =\frac12\int_0^1 \langle \nabla^*_L\bigl(v(t) \rho(t)+ \rho(t) v(t) \bigr); Y \rangle dt = \frac12\int_0^1 \langle v(t) \rho(t)+ \rho(t) v(t); \nabla_L Y \rangle dt=0.
\]
Since $Y \in {\rm ker}(\nabla_L)$ is arbitrary, we conclude that $\rho_1 -\rho_0 \in   {\rm ker}(\nabla_L)^\perp$. \endproof

\begin{remark}\label{rem:nov18.2016.8} Observe that in Lemma \ref{lem:nov18.2016.7} (ii), if we relax the assumptions on $\rho_0$ and $\rho_1$ by merely imposing that  $\rho_0, \rho_1 \in \mathcal H_+$, then   (\ref{eq:nov19.2016.2}) and (\ref{eq:nov19.2016.2b}) still imply $\rho_1 -\rho_0 \in   {\rm ker}(\nabla_L)^\perp$.
\end{remark}

For $\rho \in \cH_{++}$ and $m \in \mathbb C^{nN\times n}$ we set
\[
F(\rho, m):= {1 \over 2}\langle m,m \rho^{-1} \rangle.
\]
Given $\rho_0, \rho_1 \in  \cD_+$, denote by $\mathcal C(\rho_0, \rho_1)$ the set of paths $(\rho, v)$ such that  $\rho \in C^1([0,1], \cD_+)$ start at $\rho_0$ and end at $\rho_1,$ $v :(0,1)\rightarrow \mathcal S^N$ is Borel, $Q_{\rho}(v) \in L^1(0,1)$ and
\[
\dot \rho= {1 \over 2} \nabla_L^* (v \rho+\rho v)
\]
in the sense of distributions on $(0,1)$.  Similarly, we define $\tilde {\mathcal C}(\rho_0, \rho_1)$ to be the set of paths $(\rho, m)$ such that  $\rho \in C^1([0,1], \cD_+)$ start at $\rho_0$ and end at $\rho_1,$ $m :(0,1)\rightarrow \mathbb C^{nN \times n}$ is Borel, $F(\rho, m) \in L^1(0,1)$ and
\[
\dot \rho= {1 \over 2} \nabla_L^* (m-m_*)
\]
in the sense of distributions on $(0,1)$.

Observe that if $v \in \mathcal S^N$ and we set $m=v \rho$ then
\[
F(\rho, m)={1 \over 2} \trace(\rho v^* v).
\]
and so, the embedding $(\rho, v) \rightarrow (\rho, v \rho)$ of ${\mathcal C}(\rho_0, \rho_1)$ into $\tilde {\mathcal C}(\rho_0, \rho_1)$, extends $\frac12\trace(\rho v^* v)$ to $F(\rho, m)$. Consequently,
\[
\inf_{(\rho, v)} \Bigl\{ \int_0^1 {1 \over 2} \trace(\rho v^* v) dt \; | \; \; (\rho, v) \in \mathcal C(\rho_0, \rho_1) \Bigr\} \geq \inf_{(\rho, m)} \Bigl\{ \int_0^1 F(\rho, m) dt \; | \; \; (\rho, m) \in \tilde{\mathcal C}(\rho_0, \rho_1) \Bigr\}.
\]
We next show that the inequality can be turned into an equality.

\begin{lemma}\label{le:dec19.2016.1} If $\rho_0, \rho_1 \in  \cD_+$ then
\[
\inf_{(\rho, v)} \Bigl\{ \int_0^1 {1 \over 2} \trace(\rho v^* v)dt \; | \; \; (\rho, v) \in \mathcal C(\rho_0, \rho_1) \Bigr\} = \inf_{(\rho, m)} \Bigl\{ \int_0^1 F(\rho, m) dt \; | \; \; (\rho, m) \in \tilde{\mathcal C}(\rho_0, \rho_1) \Bigr\}
\]
\end{lemma}
\proof{} It suffices to show that for any $(\rho, m) \in \tilde{\mathcal C}(\rho_0, \rho_1)$, there exists $(\rho, v) \in \mathcal C(\rho_0, \rho_1)$ such that
\begin{equation}\label{eq:variation*3}
\int_0^1 {1 \over 2} \trace(\rho v^* v)dt \leq \int_0^1 F(\rho, m) dt.
\end{equation}
Observe that for almost every $t \in (0,1)$ we have
\begin{equation}\label{eq:variation*1}
\dot \rho(t)= {1 \over 2} \nabla_L^* (m(t)-m_*(t)) .
\end{equation}
Since both $t \rightarrow \rho(t)$ and $t \rightarrow \dot \rho(t)$ are continuous, by Corollary \ref{cor:dec18.2016.1}, there exists a continuous map $X:[0,1] \rightarrow \cH$ such that
\begin{equation}\label{eq:variation*2}
\dot \rho(t)= {1 \over 2} \nabla_L^* \Bigl(\nabla_L X(t) \rho(t)+\rho(t) \nabla_L X(t) \Bigr).
\end{equation}
Thus, for almost every $t\in (0,1)$, both (\ref{eq:variation*1}) and (\ref{eq:variation*2}) hold and so, by Proposition \ref{prop:dec18.2016.1}
\[
 {1 \over 2}\trace\bigl(\rho(t) v^*(t) v(t)\bigr) \leq F\bigl(\rho(t), m(t)\bigr)
\]
for almost every $t\in (0,1)$ with $v(t)= \nabla_L X(t).$ Thus, (\ref{eq:variation*3}) holds, which concludes the proof. \endproof

\section{Relaxation of velocity-momentum fields}\label{sec:relax}

Given $\rho_0, \rho_1 \in \cD_{+}$ we are interested in characterizing the  paths $(\rho,v)$ in  $\mathcal  H_+ \times \mathcal S^N$ that minimize the ``action integral,'' i.e., paths that possibly attain
\begin{equation}\label{eq:variation}
 \inf_{\rho \in\cH_+, v\in\cS^N} \Bigl\{ \int_0^1  \trace (\rho v^* v) dt \; \Big| \; \dot \rho ={1 \over 2} \nabla_L^*(v \rho+\rho v), \;\; \rho(0)=\rho_0,\, \rho(1)=\rho_1\Bigr\}.
\end{equation}

When $\rho>0$, Lemma \ref{le:dec19.2016.1} replaced $\trace (\rho v^* v)$ by a new expression $F(\rho, v \rho)$, introducing a new problem which, under appropriate conditions,   is a relaxation of (\ref{eq:variation}). It then becomes neccesary  to extend $F$ to the whole set $\cH \times \mathbb C^{nN \times n}$ and study the convexity properties of the extended functional. We start by introducing the open sets
\[
\mathcal O_0:=\mathcal H_{++} \times \mathbb C^{nN \times n}, \quad \mathcal O_\infty:= \{ \rho \in \mathcal H\backslash \cH_+\} \times \mathbb C^{nN \times n}.
\]
We define the functions $F, F_0, G : \mathcal H \times \mathbb C^{nN \times n}  \rightarrow [0,\infty]$ given by
\begin{equation}\label{eq:dec19.2016.6}
G(\rho, m):=\inf_{(\rho_\ell, m_\ell)}\liminf_{\ell\rightarrow \infty} \Bigl\{ {1 \over 2}\langle m_\ell; m_\ell\rho^{-1}_\ell\rangle \; | \;  (\rho_\ell, m_\ell)_\ell \subset \mathcal O_0 \; \hbox{converges to} \; (\rho, m) \Bigr\}
\end{equation}
\[
F(\rho, m)=
\left\{\begin{array}{rl}
{1 \over 2}\langle m; m\rho^{-1} \rangle & \hbox{if} \quad (\rho, m) \in \mathcal O_0\\
G(\rho, m) & \hbox{if} \quad (\rho, m) \in \mathcal H \times \mathbb C^{nN \times n} \setminus ( \mathcal O_0 \cup \mathcal O_\infty)\\
\infty & \hbox{if} \quad (\rho, m) \in \mathcal O_\infty,
\end{array}\right.
\]
and
\[
F_0(\rho, m)=
\left\{\begin{array}{rl}
{1 \over 2}\langle m; m\rho^{-1}\rangle & \hbox{if} \quad (\rho, m) \in \mathcal O_0 \\
\infty & \hbox{if} \quad (\rho, m) \in \mathcal H \times \mathbb C^{nN \times n} \setminus \mathcal O_0.
\end{array}\right.
\]

We show (cf. Lemma \ref{lem:nov29.2016.1})  that $F$ is a convex functional and then we characterize the minimizers of
\begin{equation}\label{eq:variation-new}
 \inf_{(\rho, m)} \Bigl\{ \int_0^1  F(\rho,m) dt \; \big| \; (\ref{eq:variation-neweq1}-\ref{eq:variation-neweq2}) \quad \hbox{hold} \Bigr\}.
\end{equation}
Here the infimum is performed over the set of pairs $(\rho,m)$ satisfying the requirements
\begin{equation}\label{eq:variation-neweq1}
\rho \in W^{1,2}(0,1; \cH),  \quad m \in L^{2}\bigl(0,1; \mathbb C^{nN \times n}\bigr),
\end{equation}
\begin{equation}\label{eq:variation-neweq2}
\rho(0)=\rho_0, \; \rho(1)=\rho_1, \quad \hbox{and} \quad \dot \rho ={1 \over 2} \nabla_L^*(m- m_*)
\end{equation}
in the sense of distributions on $(0,1).$

Under technical conditions, the characterizing of the minimizers $(\rho, m)$ of (\ref{eq:variation-new}) is equivalent to characterizing the minimizers  $(\rho, v)$ of (\ref{eq:variation}). We will make use of the set of paths $\lambda:[0,1] \rightarrow \cH$ such that
\begin{equation}\label{eq:variation-neweq1dual}
\lambda \in W^{1,1}(0,1; \cH),
\end{equation}
and
\begin{equation}\label{eq:variation-neweq2dual}
\dot \lambda +{1 \over 2} (\nabla_L \lambda)^*(\nabla_L \lambda) \leq 0 \quad \mbox{a.e. on} \quad (0,1).
\end{equation}

\begin{lemma}\label{lem:nov29.2016.1} The function $F$ is convex and lower semicontinuous and equals the convex envelope of $F_0.$ In addition, the Legendre transform of $F$ is
\begin{equation}\label{eq:nov29.2016.2}
F^*(a,b)=
\left\{\begin{array}{rl}
0 & \hbox{if} \quad a+ {b^*b \over 2} \leq 0\\
 & \\
\infty & \quad \mbox{otherwise}.
\end{array}\right.
\end{equation}
\end{lemma}
\proof{}
Observe that $\mathcal O_0$ is a convex set.
For $(a, b) \in \mathcal H \times \mathbb C^{nN \times n}$ we have
\[
F^*_0(a,b)= \sup_{\rho, m} \Bigl\{ \langle a; \rho \rangle +  b \cdot m - {1 \over 2}\langle m;  m\rho^{-1}\rangle \; | \; \rho>0,\; m \in \mathbb C^{nN \times n} \Bigr\}
\]
But
\[
 b \cdot  m - {1 \over 2}\langle m;  m\rho^{-1}\rangle= - {1 \over 2} \| m\rho^{-{1\over 2}} - b\rho^{1\over 2}\|^2+{1\over 2}  \langle \rho;  b^* b \rangle.
\]
Hence,
\begin{equation}\label{eq:nov29.2016.2new}
F^*_0(a,b)= \sup_{\rho} \Bigl\{ \langle a; \rho \rangle +{1\over 2}  \langle \rho; b^*b \rangle \; | \; \rho>0 \Bigr\}
=
\left\{\begin{array}{rl}
0 & \hbox{if} \quad a+ {b^*b \over 2} \leq 0\\
 & \\
\infty & \quad \hbox{otherwise}.
\end{array}\right.
\end{equation}
Denote by $F^{**}_0$ the Legendre transform of $F^*_0.$ If $(\rho, m) \in \mathcal H \times \mathbb C^{nN \times n}$ we use (\ref{eq:nov29.2016.2new}) to obtain
\begin{equation}\label{eq:nov29.2016.3}
F^{**}_0(\rho, m)=\sup_{a, b} \Bigl\{ \langle a; \rho \rangle +  b \cdot m \; | \; (a,b) \in \mathcal H \times \mathbb C^{nN \times n}, \; a+ {b^*b \over 2} \leq 0 \Bigr\}.
\end{equation}

If $(\rho, m) \in \mathcal O_0$ we can set
\[
b=m\rho^{-1}, \; a= -{1 \over 2} \rho^{-1} m^* m \rho^{-1}  \in \mathcal H.
\]
Clearly $a+ {b^*b \over 2} = 0$ and so, by (\ref{eq:nov29.2016.3})
\begin{equation}\label{eq:nov29.2016.4}
F^{**}_0(\rho, m) \geq - \langle {1 \over 2} \rho^{-1} m^* m \rho^{-1}; \rho \rangle + \langle m; m\rho^{-1}\rangle= {1 \over 2} \langle m; m\rho^{-1}\rangle= F_0(\rho, m).
\end{equation}
If $(\rho, m) \in \mathcal O_\infty$, then there exists $x \in \mathbb C^n$ such that  $\langle \rho x ; x \rangle<0.$ Set
\[ (a_\lambda,b):= (-\lambda x \otimes x, 0).\]
Observe that for any $\lambda\ge0$ we have $a_\lambda+ {b^*b\over 2} \leq 0$. Thus, by (\ref{eq:nov29.2016.3})
\begin{equation}\label{eq:nov29.2016.5}
F^{**}_0(\rho, m) \geq \lim_{\lambda \rightarrow \infty} \langle a_\lambda, \rho \rangle +  b \cdot  m = \lim_{\lambda \rightarrow \infty} -\lambda \langle \rho x ; x \rangle=\infty.
\end{equation}
Since in general $F^{**}_0 \leq F_0$, (\ref{eq:nov29.2016.4}) and (\ref{eq:nov29.2016.5}) imply that
\begin{equation}\label{eq:nov29.2016.5bbb}
F^{**}_0=F_0=F \quad \hbox{on} \quad  \mathcal O_0 \cup \mathcal O_\infty.
\end{equation}
Observe that $F$ is lower semicontinuous. We next claim that since (\ref{eq:nov29.2016.5bbb}) holds, $F$ is convex. Indeed, let $t \in (0,1)$, let $\rho^0, \rho^1 \in \cH$ and let $m^0, m^1 \in \mathbb C^{nN \times n}.$ We are to show that
\begin{equation}\label{eq:Jan.10.2017.1}
F(\rho^t, m^t) \leq (1-t)F(\rho^0, m^0)+ tF(\rho^1, m^1),
\end{equation}
where
\[
(\rho^t, m^t):=\Bigl( (1-t)\rho^0+t \rho^1, (1-t)m^0+t m^1 \Bigr).
\]
Clearly, (\ref{eq:Jan.10.2017.1}) holds if either $(\rho^0, m^0) \in \mathcal O_\infty$ or $(\rho^1, m^1) \in \mathcal O_\infty$. Since $\mathcal O_0$ is a convex set and $F^{**}_0$ is a convex function, we use  (\ref{eq:nov29.2016.5bbb}) to conclude that (\ref{eq:Jan.10.2017.1}) holds if $(\rho^0, m^0) \in \mathcal O_0$ and $(\rho^1, m^1) \in \mathcal O_0$. It remains to prove (\ref{eq:Jan.10.2017.1}) when we have either $(\rho^0, m^0) \not \in (\mathcal O_0 \cup \mathcal O_\infty)$ and $(\rho^1, m^1) \not\in \mathcal O_\infty$ or $(\rho^1, m^1) \not \in (\mathcal O_0 \cup \mathcal O_\infty)$ and $(\rho^0, m^0) \not\in \mathcal O_\infty.$ In these latter cases, there exist sequences $(\rho_\ell^0, m_\ell^0)_\ell \subset \mathcal O_0$ converging to $(\rho^0, m^0)$ and $(\rho_\ell^1, m_\ell^1)_\ell \subset \mathcal O_0$ converging to $(\rho^1, m^1)$ such that by (\ref{eq:nov29.2016.5bbb}) and the definition of $F$,
\begin{equation}\label{eq:Jan.10.2017.11/2}
F(\rho^0, m^0)= \lim_{\ell \rightarrow \infty}  F^{**}_0(\rho_\ell^0, m_\ell^0) \quad \hbox{and} \quad
F(\rho^1, m^1)= \lim_{\ell \rightarrow \infty}  F^{**}_0(\rho_\ell^1, m_\ell^1).
\end{equation}
Note that
\begin{equation}\label{eq:Jan.10.2017.2}
\bigl(\rho^t_\ell , m^t_\ell \bigr):=\Bigl((1-t) \rho_\ell^0+t\rho_\ell^1 , (1-t) m_\ell^0+t m_\ell^1\Bigr)_\ell \subset \mathcal O_0
\end{equation}
and the sequence in (\ref{eq:Jan.10.2017.2}) converges to $(\rho^t, m^t)$. Thus, using the definition of $F$,  (\ref{eq:nov29.2016.5bbb}) and the convexity property of $F^{**}_0$,  we have
\[
F(\rho^t, m^t) \leq \liminf_{\ell \rightarrow \infty} F^{**}_0\bigl(\rho^t_\ell , m^t_\ell \bigr) \leq  \liminf_{\ell \rightarrow \infty}
\left\{(1-t) F^{**}_0(\rho_\ell^0, m_\ell^0) +t F^{**}_0(\rho_\ell^1, m_\ell^1)\right\}.
\]
This, together with (\ref{eq:Jan.10.2017.11/2}) yields (\ref{eq:Jan.10.2017.1}). Thus, $F$ is convex and so, since $F$ is also lower semicontinuous, we have $F=F^{**}$.

Note that the complement of $\mathcal O_0 \cup \mathcal O_\infty$ is contained in the boundary of $\mathcal O_0 \cup \mathcal O_\infty$ and so since $F_0^{**}$ is lower semicontinuous, (\ref{eq:nov29.2016.5bbb}), in view of the definition \eqref{eq:dec19.2016.6} of $G$, implies
\begin{equation}\label{eq:nov29.2016.6}
F^{**}_0 \leq G \quad \hbox{on} \quad  \mathcal H \times \mathbb C^{nN \times n} \setminus \mathcal O_0 \cup \mathcal O_\infty.
\end{equation}
By (\ref{eq:nov29.2016.5bbb}) and (\ref{eq:nov29.2016.6}), $F_0^{**} \leq F$ and so, $F_0^{**} \leq F^{**}.$  The fact that $F \leq F_0$ yields the reversed inequality to ensure  that    $F^{**}=F_0^{**}$. As a consequence, $F^*=F_0^{***}= F_0^{*}$ and so by (\ref{eq:nov29.2016.2new}) we obtain (\ref{eq:nov29.2016.2}).
\endproof

\begin{lemma} \label{le:dec19.2016.5} Let $\rho\in\cH_+$, $X \in \cH$ and set $m:=\nabla_L X \, \rho$, then we have
\[
F(\rho, m) =   {1 \over 2} \langle \nabla_L X  \, \rho; \nabla_L X \rangle
\]
and
\begin{equation}\label{eq:dec28.2016.6--}
(a,b):= \Bigl( {1 \over 2} (\nabla_L X)^*(\nabla_L X), \nabla_L X \Bigr) \in \partial_- F(\rho, m).
\end{equation}
\end{lemma}
\proof{} For any $\epsilon>0$, we have $\rho+ \epsilon I \in \cH_{++}$ and both $\rho$ and $(\rho+ \epsilon I)^{-1}$ have the same eigenspaces. Thus $\rho$ and $(\rho+ \epsilon I)^{-1}$ commute and so, $(\rho+ \epsilon I)^{-1}\rho \in \cH$. If $\lambda_1, \cdots, \lambda_n \geq 0$ are the eigenvalues of $\rho$ then  $\lambda_1/(\lambda_1+ \epsilon), \cdots, \lambda_n /(\lambda_n+ \epsilon)\geq 0$ are those of $(\rho+ \epsilon I)^{-1}\rho$ and so, $(\rho+ \epsilon I)^{-1}\rho \in \cH_{+}.$ Thus,
\begin{equation}\label{eq:dec28.2016.6}
0 \leq \langle (\nabla_L X)^* \nabla_L X; \rho (\rho+ \epsilon I)^{-1} \rangle = \langle  \nabla_L X  (\rho+ \epsilon I)^{-1}; \nabla_L X \rho\rangle .
\end{equation}
Since $F$ is lower semicontinuous, we have
\begin{equation}\label{eq:dec28.2016.7}
 F(\rho, m) \leq \liminf_{\epsilon \rightarrow 0^+} F(\rho + \epsilon I, \nabla_L X  \, \rho) ={1 \over 2} \liminf_{\epsilon \rightarrow 0^+} \langle \nabla_L X  \, \rho ;\nabla_L X  \, \rho (\rho + \epsilon I)^{-1} \rangle
\end{equation}
But by (\ref{eq:dec28.2016.6})
\[
\langle \nabla_L X  \, \rho ;\nabla_L X  \, \rho (\rho + \epsilon I)^{-1} \rangle \leq
\langle \nabla_L X  \, \rho;\nabla_L X  \, (\rho + \epsilon I) (\rho + \epsilon I)^{-1}\rangle =  {1 \over 2} \langle \nabla_L X  \, \rho; \nabla_L X \rangle.
\]
This, together with (\ref{eq:dec28.2016.7}) implies
\begin{equation}\label{eq:dec28.2016.6-}
F(\rho, m) \leq   {1 \over 2} \langle \nabla_L X  \, \rho; \nabla_L X \rangle
\end{equation}
On the other hand, with $(a,b)$ as in \eqref{eq:dec28.2016.6--}, we have
\[
\langle a ; \rho \rangle + b\cdot m=-{1 \over 2} \langle (\nabla_L X)^* \nabla_L X ; \rho \rangle + \langle \nabla_L X ; \nabla_L X \rho \rangle= {1 \over 2}  \langle \nabla_L X \rho; \nabla_L X  \rangle
\]
We use (\ref{eq:dec28.2016.6-}) and the fact that  $F^*(a,b)=0$ (cf. by Lemma \ref{lem:nov29.2016.1}) to conclude that
\[
\langle a ; \rho \rangle + b\cdot m ={1 \over 2}  \langle \nabla_L X  \rho; \nabla_L X\rangle \geq F(\rho, m)+ F^*(a,b) \geq \langle a ; \rho \rangle + b\cdot m.
\]
Thus, $F(\rho, m)= {1 \over 2}  \langle \nabla_L X \rho; \nabla_L X \rangle$, and
\[
\langle a ; \rho \rangle +b\cdot m = F(\rho, m)+ F^*(a,b)
\]
implying \eqref{eq:dec28.2016.6--}.
\endproof

\begin{lemma}\label{rem:dec19.2016.3} We have the following:
\begin{enumerate}
\item[(i)] If $\rho \in \cH_{+}\backslash\{0\}$ and $m \in \mathbb C^{nN \times n}$ then
\begin{equation}\label{eq:dec19.2016.4}
F(\rho, m) \geq {||m||^2 \over 2\trace(\rho)}.
\end{equation}
\item[(ii)] Assume $\rho \in C([0,1]; \cH_+)$ and $m:(0,1) \rightarrow \mathbb C^{nN \times n}$ is a Borel map such that
\[
\dot \rho= {1 \over 2} \nabla_L^* (m-m_*)
\]
in the sense of distributions on $(0,1)$ and $F(\rho, m) \in L^1(0,1).$ Then $\dot \rho \in L^2(0,1; \cH)$ and there exists a constant $c_L$ independent of $(\rho, m)$ such that
\[
c_L \int_0^1 F(\rho, m) dt \geq \int_0^1 ||\dot \rho||^2 dt.
\]
Furthermore,
\[
\trace(\rho)(t) =\trace(\rho)(0).
\]
\end{enumerate}
\end{lemma}
 \proof{} (i) When $\rho\in\cH_{++}$, \eqref{eq:dec19.2016.4} is a direct consequence of the fact that  $\rho^{-1} \trace(\rho) \geq I$.
Since $F$ is defined through the liminf in \eqref{eq:dec19.2016.6}, we conclude that if $\rho \in \cH_+\backslash\{0\}$ and
$m \in \mathbb C^{nN \times n}$ then (\ref{eq:dec19.2016.4}) continue to hold.

(ii) Let $c_L$ be such that $|| \nabla_L^*(m-m_*)|| \leq 2 c_L ||m||$. Under the assumptions in (ii) we have that for almost every $t \in (0,1)$
\[
\int_0^1 ||\dot \rho||^2 dt= \int_0^1 \Bigl\|{1 \over 2} \nabla_L^* (m-m_*)\Bigr\|^2 dt \leq c_L^2 \int_0^1||m||^2dt.
\]
This, together with (i) implies the last inequality in (ii). The conservation of $\trace(\rho)(t)$ is due to the fact that $\trace(\nabla_L^*(m-m_*)) \equiv 0.$ \endproof

\section{Strong duality and conservation of the Hamiltonian}\label{sec:duality}

In this section, we state and prove our main results.We fix $\rho_0, \rho_1 \in \cD_+$ such that
\begin{equation}\label{eq:jan.10.2017.6}
\rho_0 - \rho_1 \in {\rm ker\,}(\nabla_L)^\perp.
\end{equation}
One of the aims of this section is to show under appropriate conditions that  the convex variational problems
\begin{equation}\label{eq:12.24.2016.2}
 \inf_{(\rho, m)} \Bigl\{ \int_0^1  F(\rho,m) dt \; \big| \; (\ref{eq:variation-neweq1}-\ref{eq:variation-neweq2}) \; \hbox{hold} \Bigr\}=:i_0
\end{equation}
and
\begin{equation}\label{eq:12.24.2016.2dual}
 \sup_{\lambda} \Bigl\{ \langle \lambda(1); \rho_1 \rangle- \langle \lambda(0); \rho_0 \rangle \; \big| \; (\ref{eq:variation-neweq1dual}-\ref{eq:variation-neweq2dual}) \; \hbox{hold} \Bigr\}=:j_0,
\end{equation}
are dual to each other. Recall that one of our goals is to define a Riemannian metric on $\cD_+$. In order to have finite value for $i_0$, in view of Lemma~\ref{lem:nov18.2016.7}, it is necessary to assume that $\rm{ker}(\nabla_L)$ is spanned by the identity matrix $I$. All the analysis, however, goes through without this assumption.
\begin{prop}\label{pr:dec28.2016.8} Let $\lambda$ satisfy (\ref{eq:variation-neweq1dual}-\ref{eq:variation-neweq2dual})  and   $(\rho, m)$ satisfy (\ref{eq:variation-neweq1}-\ref{eq:variation-neweq2}).
\begin{enumerate}
\item[(i)] Then
\begin{equation}\label{eq:12.24.2016.2dp1}
\langle \lambda(1); \rho_1 \rangle- \langle \lambda(0); \rho_0 \rangle \leq \int_0^1  F(\rho,m) dt.
\end{equation}
\item[(ii)] If equality holds in (\ref{eq:12.24.2016.2dp1}) then $\lambda$ is a maximizer in (\ref{eq:12.24.2016.2dual}) and $(\rho, m)$ is a minimizer in (\ref{eq:12.24.2016.2}).
\end{enumerate}
\end{prop}
\proof{} Note that (ii) is a direct consequence of (i) and so,  the only proof to supply is that of  (i). Since $\lambda$ satisfies  (\ref{eq:variation-neweq1dual}-\ref{eq:variation-neweq2dual}), we use  Lemma \ref{lem:nov29.2016.1} to conclude that $F(\dot \lambda, \nabla_L \lambda) \equiv 0$ and so,
\begin{equation}\label{eq:12.28.2016.3ba}
\langle \dot \lambda; \rho \rangle+  \nabla_L \lambda \cdot m \leq F(\rho, m)+ F^*(\dot \lambda, \nabla_L \lambda)= F(\rho, m) \quad \hbox{a.e. on} \quad (0,1).
\end{equation}
Note that by (\ref{eq:variation-neweq2})
\[
\nabla_L \lambda \cdot m= {1 \over 2} \langle \nabla_L \lambda ; m-m_* \rangle= {1 \over 2} \langle  \lambda ; \nabla_L^*( m-m_*) \rangle= \langle \lambda ; \dot \rho \rangle,
\]
 thus, (\ref{eq:12.28.2016.3ba}) implies
\begin{equation}\label{eq:12.28.2016.3baa}
{d \over dt}  \langle \lambda; \rho \rangle=  \langle \dot \lambda; \rho \rangle+  \langle \lambda ; \dot \rho \rangle \leq F(\rho, m) \quad \hbox{a.e. on} \quad (0,1).
\end{equation}
The pointwise derivative of $\langle \lambda; \rho \rangle$ coinciding with its distributional derivative, we integrate (\ref{eq:12.28.2016.3baa}) to discover that
\[
\langle \lambda(1); \rho(1) \rangle- \langle \lambda(0); \rho(0) \rangle\leq  \int_0^1 F(\rho, m)dt.
\]\endproof

By the above proposition we have the following.
\begin{remark}\label{rem:dec26.2016.1}
Assume that $\lambda \in W^{1,1}(0, 1; \cH)$ satisfies \eqref{eq:variation-neweq2dual} and $(\rho, m) \in W^{1,2}(0, 1; \cH) \times L^{2}(0, 1; \mathbb C^{nN \times n})$ is such that (\ref{eq:variation-neweq2}) holds. Then,
\[
\langle \lambda(1); \rho(1) \rangle- \langle \lambda(0); \rho(0) \rangle=\int_0^1 F(\rho,m) dt
\]
if and only if
\begin{equation} \label{eq:dec26.2016.2}
{d \over dt} \langle \lambda; \rho\rangle=F(\rho,m) \quad \mbox{a.e.}.
\end{equation}
\end{remark}

\begin{lemma}\label{le:dec22.2016.2late}
Let $(\rho,m) \in \cH_{++} \times \mathbb C^{nN \times n}.$
\begin{enumerate}
\item[(i)] The partial derivatives of $F$ with respect to $m$ and $\rho$ are
\begin{equation}\label{eq:dec22.2016.2}
\nabla_m F(\rho, m)=m \rho^{-1}\quad \hbox{and} \quad \nabla_\rho F(\rho, m)=-{1 \over 2} (m\rho^{-1})^* (m\rho^{-1}).
\end{equation}
\item[(ii)] They satisfy the relation
\begin{equation}\label{eq:dec22.2016.3bbb}
\nabla_\rho F+{1 \over 2} (\nabla_m F)^* (\nabla_m F)=0 \quad \hbox{on} \quad \mathcal O_0.
\end{equation}
\end{enumerate}
\end{lemma}
\proof{} (i) For $ r \in \cH$ with $||r||<<1$ we have
\[
(\rho+ r)^{-1}=(I+\rho^{-1} r)^{-1} \rho^{-1}= \Bigl(I +\sum_{l=1}^\infty (-1)^l (\rho^{-1} r)^l  \Bigr) \rho^{-1}= \rho^{-1}- \rho^{-1} r \rho^{-1} + o(||r||).
\]
Hence,
\[
F(\rho+r, m)={1 \over 2} \langle m ;m(\rho+ r)^{-1} \rangle= F(\rho, m)-{1 \over 2} \langle r; (m\rho^{-1})^* (m\rho^{-1}) \rangle+ o(||r||),
\]
which gives the second identity in (\ref{eq:dec22.2016.2}). The first identity is obtained in a similar manner.

(ii) By direct computations (\ref{eq:dec22.2016.3bbb}) is obtained from (\ref{eq:dec22.2016.2}).   \endproof

\begin{thm}\label{thm:dec22.2016.2} Let $\rho_0, \rho_1 \in \cD_+$ (recall that throughout this section (\ref{eq:jan.10.2017.6}) is assumed to hold).
\begin{enumerate}
\item[(i)] The  problem (\ref{eq:12.24.2016.2}) admits a minimizer $(\rho, m)$.
\item[(ii)] Let $J:=\{t \in (0,1) \; | \; {\rm det\,}(\rho(t))>0 \}.$ Then $J$ is an open set and there exists a measurable map $\hat\lambda: J \rightarrow {\rm ker}(\nabla_L)^\perp$ such that for almost every $t \in J$
\begin{equation}\label{eq:12.24.2016.3bis}
m= \nabla_L\hat\lambda\, \rho \quad \hbox{on} \quad J.
\end{equation}
\item[(iii)] If $\epsilon>0$ and we set $J_\epsilon:=\{t \in (0,1) \; | \; {\rm det\,}(\rho(t))>\epsilon \}$, then  $\hat\lambda \in L^1(J_\epsilon; \cH)$. Extend $\hat\lambda$ to $(0,1)$ by setting $\hat\lambda$ to $0$ on $(0,1)\backslash J$, and let
$\lambda(t)=\hat\lambda(t)+\Lambda(t)$ where
	\[
		\Lambda(t)=\int_0^t -\frac12 {\rm proj} \bigl((\nabla_L \hat\lambda)^* (\nabla_L \hat\lambda)\bigr)dt.
	\]
Then $m=\nabla_L\lambda \rho$ on $J$, $\lambda \in L^1(J_\epsilon; \cH)$ and
\begin{equation}\label{eq:conditionJ}
\dot \lambda +{1 \over 2} (\nabla_L \lambda)^* (\nabla_L \lambda)=0 \quad \hbox{ on} \quad J
\end{equation}
in the sense of distributions.
\end{enumerate}
\end{thm}
\proof{} (i) By Lemma \ref{lem:nov18.2016.7} and the fact that $\trace(\rho_0-\rho_1)=0$, we have $i_0<\infty.$
Let $(\rho_\ell, m_\ell)_\ell$ be a minimizing sequence of (\ref{eq:12.24.2016.2}). Using Lemma~\ref{rem:dec19.2016.3} and the fact that
\[
\sup_{\ell} \int_0^1 F(\rho_\ell, m_\ell) dt<\infty,
\]
we conclude that
\[
\trace(\rho_\ell) \equiv 1,
\]
$(m_\ell)_\ell$ is a bounded sequence in $L^2(0,1; \mathbb C^{nN \times n})$ and $(\rho_\ell)_\ell$ is a bounded sequence in
$W^{1,2}(0,1; \cH )$. Extracting subsequences if necessary, we assume without loss of generality that $(m_\ell)_\ell$ converges weakly to some $m$ in $L^2(0,1; \mathbb C^{nN \times n})$, $(\rho_\ell)_\ell$ converges strongly to some $\rho$ in $L^2(0,1; \cH)$ and $(\dot \rho_\ell)_\ell$ converges weakly to $\dot { \rho} $ in $L^2(0,1; \cH)$. Since $(\rho_\ell, m_\ell)$ satisfies (\ref{eq:variation-neweq1}-\ref{eq:variation-neweq2}), so does $( \rho, m).$ By Lemma \ref{lem:nov29.2016.1}, $F$ is convex and lower semicontinuous and so,
\[
i_0=\liminf_{\ell \rightarrow \infty} \int_0^1 F(\rho_\ell, m_\ell) dt \geq \int_0^1 F( \rho,  m) dt \geq i_0.
\]
The first equality in the above is due to the fact that $(\rho_\ell, m_\ell)_\ell$ is a minimizing sequence in (\ref{eq:12.24.2016.2}). The first inequality is due to standard results of the calculus of variations (cf. e.g. \cite{EkelandT}) which ensure lower semicontinuity of functionals for weak topologies. The last inequality is due to the definition of $i_0$.
This proves that $(\rho, m)$ is a minimizer in (\ref{eq:12.24.2016.2}).

(ii) Since $ \rho \in W^{1,2}(0, 1; \cH_+)$,  $t \rightarrow {\rm det\,}(\rho(t))$ is a continuous function on $[0,1]$ and so, the set $J$ is an open set. The last identity in (\ref{eq:variation-neweq2}), which holds in the sense of distributions, also holds pointwise almost everywhere. Hence, for almost every $t \in J$, $m(t)$ minimizes $F(\rho(t), w)$ over the set of $w \in \mathbb C^{nN \times n}$ such that  such that
\[
\dot \rho(t)=\frac12\nabla_L^*(w-w_*).
\]
By Proposition \ref{prop:dec18.2016.1}, for these $t,$ there exists a unique $\hat\lambda(t) \in {\rm ker}(\nabla_L)^\perp$ such that $m(t)=\nabla_L \hat\lambda(t) \rho(t).$ By Corollary \ref{cor:dec18.2016.1}, the map $\hat\lambda:J \rightarrow {\rm ker}(\nabla_L)^\perp$ is measurable.

(iii) We first establish $\hat\lambda \in L^1(J_\epsilon, \cH)$. By Lemma~\ref{rem:dec19.2016.3} (i) and the fact that $F(\rho, m) \in L^1(0,1)$, we have $m \in L^2(0,1; \mathbb C^{nN \times n})$, and therefore $\nabla_L \hat\lambda \in  L^2(J_\epsilon; \mathbb C^{nN \times n}).$ It follows
\[
||\nabla_L \hat\lambda||^2 \in  L^1(J_\epsilon; \mathbb C^{nN \times n}).
\]
Now, we apply the Poincar\'e--Wirtinger inequality (cf. Theorem \ref{thm:nov18.2016.3}) with $\mathbb K:=\{I\}$ to conclude that there exists a constant $c_{\mathbb K}$ independent of $\epsilon$ and $\hat\lambda$ such that
\[
c_{\mathbb K} \int_{J_\epsilon}||\hat\lambda||^2 dt \leq  \int_{J_\epsilon} ||\nabla_L \hat\lambda||^2 dt.
\]
Therefore, $\hat\lambda \in L^1(J_\epsilon; \cH)$. Since $\nabla_L \hat\lambda \in  L^2(J_\epsilon; \mathbb C^{nN \times n})$, we have $\Lambda \in L^1(0,1)$. It follows $\lambda \in L^1(J_\epsilon; \cH)$. Recalling that $\Lambda(t) \in {\rm ker}(\nabla_L)$ for any $t\in (0,1)$, we have $m=\nabla_L \hat\lambda \rho=\nabla_L\lambda \rho$.
Proving \eqref{eq:conditionJ} amounts to proving that, for any arbitrary $f \in C_c^1(J; \cH)$,
\begin{equation}\label{eq:12.24.2016.3ter}
\int_J \langle \dot f; \lambda\rangle dt= {1 \over 2} \int_J \langle f ;(\nabla_L \lambda)^* (\nabla_L \lambda) \rangle dt.
\end{equation}
For $f \in C_c^1(J; {\rm ker}(\nabla_L))$, since $\hat\lambda(t) \in {\rm ker}(\nabla_L)^\perp$, $\Lambda(t) \in {\rm ker}(\nabla_L)$ and $ f(t) \in {\rm ker}(\nabla_L)$,  we have
\[
\int_J \langle \dot f; \lambda\rangle dt= \int_J \langle \dot f; \Lambda\rangle dt =\int_J \langle f; {\rm proj}\Bigl({1 \over 2} (\nabla_L \hat \lambda)^* (\nabla_L \hat \lambda)\Bigr)\rangle dt=
\int_J \langle f; \Bigl({1 \over 2} (\nabla_L \hat \lambda)^* (\nabla_L \hat \lambda)\Bigr)\rangle dt ,
\]
which proves (\ref{eq:12.24.2016.3ter}). Therefore, it is remains to consider $f \in C_c^1(J; {\rm ker}(\nabla_L)^\perp)$.
Fix such an $f$ and   denote the support of $f$ by ${\rm spt}(f).$ To avoid technical difficulties, we assume without loss of generality that ${\rm spt}(f)$ is contained in some $J_\epsilon.$ Extend $f$ to $[0,1]$ by setting $f(t)\equiv 0$ on $[0,1] \setminus J$ and observe that the extension, which we still denote by $f,$ satisfies $f \in C_c^1([0,1]; {\rm ker}(\nabla_L)^\perp)$. By Proposition \ref{prop:dec18.2016.1} and Corollary \ref{cor:dec18.2016.1}, there exists a unique map $\beta \in C(J ; {\rm ker}(\nabla_L)^\perp)$ such that
\[
\dot f= {1 \over 2} \nabla_L^*( \nabla_L \beta \rho+ \rho \nabla_L \beta) \quad \hbox{on} \quad J.
\]
By its uniqueness property on $J$, we have $\beta(t)=0$ for $t \in J \setminus {\rm spt}(f).$ Set $\beta(t)=0$ for $t \in [0,1] \setminus J$ and continue to denote the extension by $\beta$ to observe that  $\beta \in C([0,1] ; {\rm ker}(\nabla_L)^\perp)$ and
\begin{equation}\label{eq:12.24.2016.3ter2}
\dot f= {1 \over 2} \nabla_L^*( \nabla_L \beta \rho+ \rho \nabla_L \beta) \quad \hbox{on} \quad (0,1).
\end{equation}
We set
\[
\rho_\epsilon:= \rho+ \epsilon f, \quad m_\epsilon:=m+\epsilon \nabla_L \beta \rho.
\]
Since ${\rm spt}(f)$ is a compact subset of $J$ there exists $c>0$ such that $\rho \geq c$ on ${\rm spt}(f)$. We have
\[
0 \leq \int_0^1 F(\rho_\epsilon, m_\epsilon)dt -\int_0^1 F(\rho, m)dt  =\int_{{\rm spt}(f)} \bigl( F(\rho_\epsilon, m_\epsilon) - F( \rho, m)\bigr)dt
\]
and so, the function $\epsilon \rightarrow \int_0^1 F(\rho_\epsilon, m_\epsilon)dt$ achieves its minimum at $\epsilon=0$. Since $ m \in L^2(0, 1; \mathbb C^{nN \times n})$ and $F$ is differentiable on $\{r \in \cH\; | \; r \geq c \} \times \mathbb C^{nN \times n}$ with its derivatives given by (\ref{eq:dec22.2016.2}), we conclude that $\int_0^1 F(\rho_\epsilon, m_\epsilon)dt -\int_0^1 F( \rho, \ m)dt$ is differentiable  at $\epsilon=0$ with a null derivative there. More precisely,
\begin{eqnarray}
0 &=&
\int_{{\rm spt}(f)} \bigl( \langle \nabla_\rho F( \rho, \ m); f \rangle +  \nabla_m F( \rho, m)\cdot \nabla_L \beta \rho\bigr)dt \\
&=&  \int_{J} \bigl( \langle \nabla_\rho F(\rho, m); f \rangle +  \nabla_m F( \rho, m)\cdot \nabla_L \beta \rho \bigr)dt .
\end{eqnarray}
This, together with (\ref{eq:dec22.2016.2}) and the fact that $m= \nabla_L \lambda \rho$ on $J$ yields
\begin{eqnarray*}
0 &=&  \int_{J} \bigl( \langle -{1 \over 2} (\nabla_L \lambda )^* \nabla_L \lambda; f \rangle + \nabla_L \lambda \cdot \nabla_L \beta \rho \bigr)dt \\
&=& \int_{J} \bigl( \langle - {1 \over 2} (\nabla_L \lambda )^* \nabla_L \lambda; f \rangle + {1 \over 2} \langle \nabla_L \lambda ; \nabla_L \beta \rho + \rho \nabla_L \beta  \rangle \bigr)dt\\
&=& \int_{J} \bigl( \langle - {1 \over 2} (\nabla_L \lambda )^* \nabla_L \lambda; f \rangle + {1 \over 2} \langle \lambda ; \nabla_L^* \bigl(\nabla_L \beta \rho + \rho \nabla_L \beta \bigr) \rangle \bigr)dt.
\end{eqnarray*}
We then use (\ref{eq:12.24.2016.3ter2}) to obtain (\ref{eq:12.24.2016.3ter}). \endproof

\begin{cor} \label{cor:dec28.2016.1}  Let $\rho_0, \rho_1 \in \cD_+$ and let $(\rho, m)$  be such that (\ref{eq:variation-neweq1}-\ref{eq:variation-neweq2})holds.
\begin{enumerate}
\item[(i)] If there exists $\lambda \in W^{1,1}(0,1; \cH)$ such that
\begin{equation}\label{eq:12.28.2016.1}
\dot \lambda +{1 \over 2} (\nabla_L \lambda)^* (\nabla_L \lambda)=0 \quad \hbox{on} \quad (0,1)
\end{equation}
in the sense of distributions and
\begin{equation}\label{eq:12.28.2016.1new}
m= \nabla_L \lambda \rho \quad \hbox{a.e. on} \quad (0,1),
\end{equation}
 then $(\rho, m)$ minimizes (\ref{eq:12.24.2016.2}).
\item[(ii)] Conversely, assume that $(\rho, m)$ minimizes (\ref{eq:12.24.2016.2}) and the range of $\rho$ is contained in $\cD_+$. Then, there exists $\lambda \in W^{1,1}(0,1; \cH)$ such that (\ref{eq:12.28.2016.1})  holds.
\item[(iii)] Any minimizer $(\rho, m)$ of (\ref{eq:12.24.2016.2}) whose range is contained in $\cD_+$ must be of class $C^\infty.$
\end{enumerate}
\end{cor}
\proof{} (i) Assume there exists $\lambda \in W^{1,1}(0,1; \cH)$ such that (\ref{eq:12.28.2016.1}) holds.  Since (\ref{eq:12.28.2016.1}) holds almost everywhere and $m$ satisfies (\ref{eq:12.28.2016.1new}), we have
\[
\langle \dot \lambda; \rho \rangle+ \nabla_L \lambda \cdot m=-\langle {1 \over 2} (\nabla_L \lambda)^* (\nabla_L \lambda); \rho \rangle+
\langle \nabla_L \lambda ;   \nabla_L \lambda \rho \rangle   = {1 \over 2}  \langle \nabla_L \lambda ;   \nabla_L \lambda \rho \rangle
\]
Hence by Lemma \ref{le:dec19.2016.5} we have
\begin{equation}\label{eq:12.28.2016.3}
\langle \dot \lambda; \rho \rangle+  \nabla_L \lambda \cdot m = F(\rho, m)+ F^*(\dot \lambda, \nabla_L \lambda)= F(\rho, m) \quad \hbox{a.e. on} \quad (0,1)
\end{equation}
Since
\[
\nabla_L \lambda \cdot m = {1 \over 2} \langle \nabla_L \lambda ; m-m_* \rangle= {1 \over 2} \langle  \lambda ; \nabla_L^*( m-m_*) \rangle,
\]
we combine (\ref{eq:variation-neweq2}) and (\ref{eq:12.28.2016.3}) to conclude that
\begin{equation}\label{eq:12.28.2016.3b}
{d \over dt}  \langle \lambda; \rho \rangle=  \langle \dot \lambda; \rho \rangle+  \langle \lambda ; \dot \rho \rangle= F(\rho, m)+ F^*(\dot \lambda, \nabla_L \lambda)= F(\rho, m) \quad \hbox{a.e. on} \quad (0,1).
\end{equation}
The pointwise derivative of $\langle \lambda;\rho \rangle$ coinciding with its distributional derivative, we integrate (\ref{eq:12.28.2016.3b}) to discover that
\begin{equation}\label{eq:12.28.2016.3c}
\langle \lambda(1); \rho(1) \rangle- \langle \lambda(0); \rho(0) \rangle=  \int_0^1 F(\rho, m)dt.
\end{equation}
We use Proposition \ref{pr:dec28.2016.8} to conclude that $(\rho, m)$ minimizes (\ref{eq:12.24.2016.2}).

(ii)  Assume that $(\rho, m)$ minimizes (\ref{eq:12.24.2016.2}) and the range of $\rho$ is contained in $\cD_+$.  By Theorem \ref{thm:dec22.2016.2}, there exists $\lambda: [0,1] \rightarrow \cH$ such that (\ref{eq:12.28.2016.1new}) holds. Since $\rho$ is continuous, its range is a compact set and so, the range of ${\rm det\,}(\rho)$ is a compact subset of $(0,\infty)$. Since $m \in L^2(0,1; \mathbb C^{nN \times n})$ we have $\nabla_L \lambda \in L^2(0,1; \mathbb C^{nN \times n})$. Thus, in view of \eqref{eq:12.28.2016.1}, $\lambda \in W^{1,1}(0,1; \cH)$.

(iii) Assume $(\rho, m)$ is a minimizer in (\ref{eq:12.24.2016.2}) and the range of $\rho$ is contained in $\cD_+$ . Since by (ii) $\lambda$ is continuous, (\ref{eq:12.28.2016.1}) implies that $\dot \lambda$ is continuous and so, $\lambda$ is of class $C^1.$ We repeat the procedure to conclude that $\lambda$ is of class $C^\infty.$ Since (\ref{eq:12.28.2016.1new}) holds and both $\rho$ and $\lambda$ are continuous, we obtain that $m$ is continuous. By (\ref{eq:variation-neweq2}), $\dot \rho$ is continuous and so, $\rho$ is of class $C^1.$ Because, $\lambda$ has been shown to be of class $C^\infty$, (\ref{eq:12.28.2016.1new}) implies that $m$ is of class $C^1.$ We use again (\ref{eq:variation-neweq2}) to conclude that $\dot \rho$ is of class $C^1$ and so, $\rho$ is of class $C^2.$ We repeat the procedure to conclude that  $\rho$ is of class $C^\infty.$ \endproof

\begin{remark}\label{rem:dec28.2016.10} Let $\rho_0, \rho_1 \in \cD_+$. By Theorem \ref{thm:dec22.2016.2},    (\ref{eq:12.24.2016.2}) admits a minimizer $(\tilde \rho, \tilde m)$. Observe that thanks to Corollary \ref{cor:dec28.2016.1}, we have proven that if the range of $\tilde \rho$ is contained in $\cH_{++}$, then we have the duality result
\[
 \min_{(\rho, m)} \Bigl\{ \int_0^1  F(\rho,m) dt \; \big| \; (\ref{eq:variation-neweq1}-\ref{eq:variation-neweq2}) \; \hbox{hold} \Bigr\}=
 \max_{\lambda} \Bigl\{ \langle \lambda(1); \rho_1 \rangle- \langle \lambda(0); \rho_0 \rangle \; \big| \; (\ref{eq:variation-neweq1dual}-\ref{eq:variation-neweq2dual}) \; \hbox{hold}
 \Bigr\}.
 \]
\end{remark}

Our goal is to extend the duality result in Remark \ref{rem:dec28.2016.10} without assuming that the range of $\tilde \rho$ is contained in $\cH_{++}$, at some expense.  It is convenient to introduce the sets
\[
\mathcal A:=  \{ \rho \in L^2(0,1; \cH) \; | \; \trace(\rho) \equiv 1\} \times L^2(0,1; \mathbb C^{nN \times n}),
\]
\[
\mathcal A_1:=  \{ \rho \in L^2(0,1; \cH) \; | \; \trace(\rho) \leq 1\} \times L^2(0,1; \mathbb C^{nN \times n}),
\]
\[
\mathcal A_\infty:= L^2(0,1; \cH) \times L^2(0,1; \mathbb C^{nN \times n}),
\]
and
\[
\mathcal B:= W^{1,2}(0,1; \mathcal H), \quad \mathcal B_\ell:= \Bigl\{ \lambda \in W^{1,2}(0,1; \mathcal H)\; | \; ||\lambda||_{W^{1,2}} \leq \ell^2 \Bigr\},
\]
where
\[
||\lambda||_{W^{1,2}}^2:= \int_0^1 (||\lambda||^2+||\dot \lambda||^2)dt .
\]
We also set
\[
J(a, b)=
\inf_{(\rho,m)} \bigl\{  F(\rho, m)- \langle \rho;a \rangle-{1 \over 2} \langle m-m_* ;b \rangle  \; | \; (\rho,m) \in \mathcal A \bigr\}.
\]
and for $\beta \in \{1, \infty\}$,
\[
J_\beta(a, b)=
\inf_{(\rho,m)} \bigl\{  F(\rho, m)- \langle \rho;a \rangle-{1 \over 2} \langle m-m_* ;b \rangle  \; | \; (\rho,m) \in \mathcal A_\beta \bigr\}.
\]

\begin{remark}\label{rem:june18.2017.2new} Let $\lambda \in W^{1,2}(0,1; \mathcal H)$, let  $\alpha \in W^{1,2}(0,1)$ and set $\bar \lambda:=\lambda+ \alpha I,$
where $I$ is the identity matrix. Then
\begin{enumerate}
\item[(i)] Since $F$ is $1$--homogeneous, $J_\infty=-F^* \quad \hbox{and} \quad J_1=-\sup_{0\le \mu \le 1} \{-\mu J\}=-J_{-}.$
\item[(ii)] $J(\dot \lambda, \nabla_L \lambda) \in L^2(0,1)$.
\item[(iii)]  Since $I \in \rm{ker}(\nabla_L),$ $\nabla_L \lambda= \nabla_L \bar \lambda.$ One may easily check that $J(\dot {\bar \lambda}, \nabla_L \bar \lambda)=J(\dot \lambda, \nabla_L \lambda) -\dot \alpha.$
\end{enumerate}
\end{remark}

\begin{lemma}\label{lem:june18.2017.6new} For any $\lambda \in W^{1,2}(0,1; \mathcal H)$, there exists $\bar \lambda \in W^{1,2}(0,1; \mathcal H)$  such that $J(\dot {\bar \lambda}, \nabla_L \bar \lambda) \geq 0$ and
\[
\inf_{(\rho, m) \in \mathcal A_1} \mathcal L(\rho,m, \lambda)= \langle \bar \lambda(1); \rho_1\rangle- \langle \bar \lambda(0); \rho_0\rangle.
\]
\end{lemma}
\proof{} Set
\[
F:=\{t \in (0,1) \; |\; J(\dot \lambda(t), \nabla_L \lambda(t))<0  \}, \quad \alpha(t):= \int_0^t \chi_F(s) J(\dot \lambda(s), \nabla_L \lambda(s)) ds \quad \forall \; t \in (0,1).
\]
By Remark \ref{rem:june18.2017.2new}, $\alpha \in W^{1,2}(0,1)$ and $\bar \lambda:= \lambda + \alpha I$ satisfies the desired properties. \endproof
\begin{cor}\label{cor:june18.2017.8} We have
\[
\sup_{\lambda \in \mathcal B} \inf_{(\rho, m) \in \mathcal A} \mathcal L(\rho,m, \lambda)=\sup_{\lambda \in \mathcal B} \inf_{(\rho, m) \in \mathcal A_1} \mathcal L(\rho,m, \lambda)= \sup_{\lambda \in \mathcal B} \inf_{(\rho, m) \in \mathcal A_\infty} \mathcal L(\rho,m, \lambda)
\]
\end{cor}
\proof{} Since by Lemma \ref{lem:june18.2017.6new}
\[
\sup_{\lambda \in \mathcal B} \inf_{(\rho, m) \in \mathcal A_1} \mathcal L(\rho,m, \lambda) =
\sup_{\lambda \in \mathcal B} \Bigl\{ \inf_{(\rho, m) \in \mathcal A_1} \mathcal L(\rho,m, \lambda)\; | \; J(\dot \lambda, \nabla_L \lambda) \geq 0\Bigr\}
\]
we use Remark \ref{rem:june18.2017.2new} to conclude that
\[
\sup_{\lambda \in \mathcal B} \inf_{(\rho, m) \in \mathcal A_1} \mathcal L(\rho,m, \lambda) =
\sup_{\lambda \in \mathcal B} \Bigl\{ \langle \lambda(1); \rho_1\rangle- \langle  \lambda(0); \rho_0\rangle\; | \; J(\dot \lambda, \nabla_L \lambda) \geq 0\Bigr\}.\]
Similarly, Lemma~\ref{lem:june18.2017.6new} and Remark~\ref{rem:june18.2017.2new} imply that
\[ \inf_{(\rho, m) \in \mathcal A_1} \mathcal L(\rho,m, \lambda)= \inf_{(\rho, m) \in \mathcal A} \mathcal L(\rho,m, \lambda).\]
\endproof

\begin{thm}\label{thm:dec22.2016.2newnew} Let $\rho_0, \rho_1 \in \cD_+$. We have
\[
 \min_{(\rho, m) \in \mathcal A} \Bigl\{ \int_0^1  F(\rho,m) dt \; \big| \; \eqref{eq:variation-neweq2} \; \hbox{holds} \Bigr\}=
 \sup_{\lambda \in \mathcal B} \Bigl\{ \langle \lambda(1); \rho_1 \rangle- \langle \lambda(0); \rho_0 \rangle \; \big| \;\eqref{eq:variation-neweq2dual} \; \hbox{holds} \Bigr\}.
\]
\end{thm}
\proof{} We endow $\mathcal A$ and $\mathcal B$ with their respective weak topologies and for $(\rho, m) \in \mathcal A$ and $\lambda \in \mathcal B$ we define
\[
\mathcal L(\rho,m, \lambda):= \langle \lambda(1); \rho_1\rangle- \langle \lambda(0); \rho_0\rangle +\int_0^1 \bigl(F(\rho, m)- \langle \rho;\dot \lambda \rangle-{1 \over 2} \langle m-m_* ;\nabla_L \lambda \rangle
\bigr)dt.
\]
For $\ell \in (0,\infty)$, $\mathcal B_\ell$ is a compact convex topological space. Let $(\rho^0, m^0) \in \mathcal A$ and $\lambda^0 \in \mathcal B_\ell.$   For any $c$, the set $\{\lambda \in \mathcal B_\ell \; | \;  \mathcal L(\rho^0,m^0, \lambda) \geq c\}$ a closed convex set in $\mathcal B_\ell $ while the set $\{(\rho, m) \in \mathcal A \; | \;  \mathcal L(\rho,m, \lambda^0) \leq c\}$ a closed convex set in $\mathcal A $. Thus, by Theorem 1.6 \cite{Mertens}
\begin{equation}\label{eq:dec31.2016.5}
\inf_{\mathcal A}\sup_{\mathcal B_\ell} \mathcal L= \sup_{\mathcal B_\ell} \inf_{\mathcal A} \mathcal L
\end{equation}
Set
\[
\mathcal E(\rho, m):= \sup_{\lambda \in \mathcal B_1} \langle \lambda(1); \rho_1\rangle- \langle \lambda(0); \rho_0\rangle - \int_0^1 \bigl( \langle \rho;\dot \lambda \rangle+{1 \over 2} \langle m-m_* ;\nabla_L \lambda \rangle)dt.
\]
Then $\mathcal E$ is a nonnegative convex function such that
\begin{equation}\label{eq:dec31.2016.6}
\mathcal E(\rho, m)=0 \quad \iff \quad  \rho(1)=\rho_1, \; \rho(0)=\rho_0 \; \hbox{and} \;\; \dot \rho ={1 \over 2} \nabla_L^*(m-m_*)
\end{equation}
in the sense of distributions on $(0,1).$

For any $(\rho, m) \in \mathcal A$
\[
\sup_{\lambda \in \mathcal B_\ell} \mathcal L(\rho, m)=\int_0^1 F(\rho, m) dt + \ell \mathcal E(\rho, m)
\]
Let $(\rho_\ell, m_\ell) \in \mathcal A$ be such that
\[
\inf_{(\rho,m) \in \mathcal A}\sup_{\lambda \in \mathcal B_\ell} \mathcal L(\rho, m, \lambda)=  \int_0^1 F(\rho_\ell, m_\ell) dt + \ell \mathcal E(\rho_\ell, m_\ell)
\]
Since $\trace(\rho_\ell) \equiv 1$, by Lemma~\ref{rem:dec19.2016.3} $(m_\ell)_\ell$ is bounded in $L^2(0,1; \mathbb C^{nN \times n})$. The fact that $\rho_\ell \geq 0$ yields that $(\rho_\ell)_\ell$ is bounded in $L^2(0,1; \cH).$ Thus, there exists a subsequence $(  \rho_{\ell_k}, m_{\ell_k} )_k$ such that $(\rho_{\ell_k})_k$ converges weakly to some $\rho_\infty$ is $L^2(0,1; \cH)$ and $(m_{\ell_k})_k$ converges weakly to some $m_\infty$ is $L^2(0,1;  \mathbb C^{nN \times n}).$ Clearly, we have $\trace(\rho_\infty) \equiv 1.$

Let $(\tilde \rho, \tilde m)$ be a minimizer of (\ref{eq:12.24.2016.2}) as given by Theorem \ref{thm:dec22.2016.2newnew}. By (\ref{eq:dec31.2016.6}), $\mathcal E(\tilde\rho, \tilde m)=0$ and so,
\begin{equation}\label{eq:dec31.2016.7}
\int_0^1 F(\rho_\ell, m_\ell) dt + \ell \mathcal E(\rho_\ell, m_\ell) \leq \int_0^1 F(\tilde \rho, \tilde m) dt.
\end{equation}
Hence, by the weak lower semicontinuity property of $\mathcal E$ we have
\[
\mathcal E(\rho_\infty, m_\infty) \leq \liminf_{k \rightarrow \infty } \mathcal E(\rho_{\ell_k}, m_{\ell_k}) \leq  \liminf_{k \rightarrow \infty } {1 \over \ell_k}  \int_0^1 F(\tilde \rho, \tilde m) dt=0
\]
We conclude that $\mathcal E(\rho_\infty, m_\infty)=0$ and so by (\ref{eq:dec31.2016.6})
\begin{equation}\label{eq:dec31.2016.8}
\rho_\infty(1)=\rho_1, \; \rho_\infty(0)=\rho_0 \; \hbox{and} \;\;  \dot \rho_\infty ={1 \over 2} \nabla_L^*(m_\infty-(m_\infty)_*)
\end{equation}
in the sense of distributions on $(0,1).$ By (\ref{eq:dec31.2016.7})
\begin{equation}\label{eq:dec31.2016.9}
\int_0^1 F(\rho_\infty, m_\infty) dt \leq \liminf_{k \rightarrow \infty } \int_0^1 F(\rho_{\ell_k}, m_{\ell_k}) dt \leq \int_0^1 F(\tilde \rho, \tilde m) dt.
\end{equation}
Since $(\rho_\infty, m_\infty)$ satisfies (\ref{eq:dec31.2016.8}), (\ref{eq:dec31.2016.9}) shows that its is also a minimizer in (\ref{eq:12.24.2016.2}).

By  the definition of $(\rho_\ell, m_\ell)$ and then (\ref{eq:dec31.2016.5}), we have
\[
\int_0^1 F(\rho_{\ell_k}, m_{\ell_k}) dt + \ell_k \mathcal E(\rho_{\ell_k}, m_{\ell_k}) = \sup_{\mathcal B_{\ell_k}} \inf_{\mathcal A} \mathcal L \leq  \sup_{\mathcal B} \inf_{\mathcal A} \mathcal L
\]
and so,
\[
\int_0^1 F(\rho_{\ell_k}, m_{\ell_k}) dt  \leq  \sup_{\mathcal B} \inf_{\mathcal A} \mathcal L.
\]
This, together with (\ref{eq:dec31.2016.9}) and Corollary~\ref{cor:june18.2017.8}, implies
\[
\int_0^1 F(\rho_\infty, m_\infty) dt \leq \sup_{\mathcal \lambda \in B} \inf_{ (\rho, m) \in \mathcal A} \mathcal L= \sup_{\lambda \in \mathcal B} \inf_{(\rho, m) \in \mathcal A_\infty} \mathcal L(\rho,m, \lambda)=\sup_{\lambda \in \mathcal B}  \Bigl\{\langle \lambda(1); \rho_1\rangle- \langle \lambda(0); \rho_0\rangle- \int_0^1 F^*(\dot \lambda, \nabla_L \lambda) dt\Bigr\}.
\]
Hence, using Lemma~\ref{lem:nov29.2016.1} we conclude that
\[
\int_0^1 F(\rho_\infty, m_\infty) dt \leq  \sup_{\lambda \in \mathcal B}  \Bigl\{\langle \lambda(1); \rho_1\rangle- \langle \lambda(0); \rho_0\rangle\; | \; \dot \lambda +{1 \over 2} (\nabla_L \lambda)^* (\nabla_L \lambda) \leq 0\Bigr\}.
\]
This, together with Proposition \ref{pr:dec28.2016.8} yields the desired result.
\endproof

\begin{thm}[Conservation of the Hamiltonian]\label{lem:nov20.2016.7} Let  $\rho_0, \rho_1 \in \cD_+$ and assume $(\rho, m)$ minimizes  (\ref{eq:12.24.2016.2}). Then

\noindent (i)
\[
F ( \rho(t),  m(t)) \equiv F ( \rho(0),  m(0)).
\]

\noindent (ii)  If $0 \leq s \leq t \leq 1$ then
\[
W_{2}(\rho(s), \rho(t))=(t-s) \sqrt{2F( \rho(t), m(t))}=(t-s) W_{2}(\rho_0, \rho_1).
\]

\noindent (iii) If we further assume that $\lambda \in W^{1,1}(0,1; \cH)$ is a maximizer in (\ref{eq:12.24.2016.2dual}) then
\[
\langle \lambda(t) ; \rho(t) \rangle= \langle \lambda(0) ; \rho_0 \rangle +   {W_{2}(\rho_0, \rho(t))^2 \over 2t}, ~~t\in (0,1].
\]
\end{thm}
\proof{} (i) Let $\zeta \in C_c^1(0,1)$ be arbitrary and set $S(t)=t+\epsilon \zeta(t).$ We have $S(0)=0,$ $S(1)=1$ and $\dot S(t)=1+\epsilon \dot \zeta(t)>1/2$ for $|\epsilon|<<1.$ Thus, $S:[0,1] \rightarrow [0,1]$ is a diffeomorphism. Let $T:=S^{-1}$ and set
\[
f(s)= \rho(T(s)), \quad w(s)= \dot T(s) m(T(s)).
\]
We have
\[
\dot f ={1 \over 2} \nabla_L^* (w-w_*), \quad f(0)=\rho_0, \quad f(1)=\rho_1.
\]
Thus,
\[
\int_0^1 F( \rho, m) dt \leq \int_0^1 F ( f, w) ds=  \int_0^1 \dot T^2(s) F( \rho(T(s)), m(T(s))) ds.
\]
We use the fact that $dt= \dot T(s) ds$ and $\dot T(S(t)) \dot S(t)=1$ to conclude that
\[
\int_0^1 F( \rho, m) dt  \leq \int_0^1 {1 \over \dot S(t)} F ( \rho(t), m(t)) dt= \int_0^1 (1- \epsilon \dot \zeta(t) + o(\epsilon)) F( \rho(t), m(t)) dt.
\]
Since $\epsilon \rightarrow \int_0^1 (1- \epsilon \dot \zeta(t) + o(\epsilon)) F ( \rho(t), m(t)  dt$ admits its minimum at $0$, we conclude that its derivative there is null, i.e.,
\[
\int_0^1 \dot \zeta(t) F ( \rho(t),  m(t)) dt=0.
\]
This proves that the distributional derivative of $F ( \rho(t),  m(t))$ is null and so, $F ( \rho(t),  m(t))$ is independent of $t$.

(ii) Recalling the definition of $W_2$ in \eqref{eq:quantumomt}, we have
	\[
		W_2(\rho_0,\rho_1)^2=\int_0^1 2F(\rho(\tau),m(\tau))d\tau.
	\]
Due to the time homogeneity of the definition, and the optimality of $(\rho, m)$, one can clearly see, for $0 \leq s \leq t \leq 1$,
	\[
		W_2(\rho(s),\rho(t))^2=(t-s)\int_s^t 2F(\rho(\tau),m(\tau)) d\tau.
	\]
We use these, together with (i), to conclude the proof of (ii).

(iii) We use Remark \ref{rem:dec26.2016.1} and the duality result in Theorem \ref{thm:dec22.2016.2newnew}  to conclude that
\[
{d \over dt} \langle \rho; \lambda \rangle= F(\rho, m) \quad \mbox{a.e.}.
\]
Thus, if $0\leq s<t \leq 1$, then
\[
\langle \lambda(t); \rho(t) \rangle- \langle \lambda(s); \rho(s) \rangle= \int_s^t F(\rho(\tau), m(\tau)) d\tau.
\]
We apply (i) and use (ii) to conclude that
\[
\langle \lambda(t); \rho(t) \rangle- \langle \lambda(s); \rho(s) \rangle= (t-s) F ( \rho(0),  m(0))= {W_{2}(\rho(s), \rho(t))^2 \over 2(t-s)}.
\] \endproof

\section{Conclusions and further research}\label{sec:conclusions}

This note continues our study of a quantum mechanical approach to (non-commutative) optimal mass transport between density matrices initiated in \cite{Yongxin}. In particular, we prove a duality result that elucidates the connection of our set-up to Monge-Kantorovich theory \cite{Kantorovich1948}, in particular Kantorovich duality as well as a Poincare\'e-Wirtinger type result. For applications, it is important to note that our methodology leads to convex optimization problems that may be implemented and numerically solved on computer.

It is of interest to explore further potential implications of this construction to quantum channels and quantum information. It would seem that our results would seem to rule out the possibility of naturally defining a joint probabilistic structure across a (possibly unknown) quantum channel with known marginal density matrices, and the best one could do is along the lines of Theorem~\ref{lem:nov20.2016.7}. In this sense, it may be that the dynamic Benamou-Brenier approach to mass transport may be the more versatile formulation of defining the  Wasserstein metric than the classical Monge-Kantorovich approach. Finally, much of the theory should go through in the infinite dimensional case. This is another area we plan to further explore.

\section*{Acknowledgements}

This project was supported by AFOSR grants FA9550-15-1-0045 and FA9550-17-1-0435, grants from the National Center for Research Resources P41-RR-013218 and the National Institute of Biomedical Imaging and Bioengineering P41-EB-015902, National Science Foundation Grant DMS-1160939, grants from the National Institutes of Health 1U24CA18092401A1 and R01-AG048769, and a postdoctoral fellowship through Memorial Sloan Kettering Cancer Center.

\bibliographystyle{plain}

\end{document}